\documentclass[a4paper,11pt]{amsart}
\usepackage{amsmath}
\usepackage{amssymb}
\usepackage{amsfonts}

\def \endproof {\quad \hfill  \rule{2mm}{2mm} \par\medskip}

\DeclareMathSymbol{\subsetneqq}{\mathbin}{AMSb}{36}

\def \endproof {\quad \hfill  \rule{2mm}{2mm} \par\medskip}
\DeclareMathSymbol{\subsetneqq}{\mathbin}{AMSb}{36}

\newcommand{\R}{\mathbb{R}}
\newcommand{\Z}{\mathbb{Z}}

\newtheorem{th1}{{\bf Theorem}}[section]

\newtheorem{lem}{{\bf Lemma}}[section]
\newtheorem{prop}{{\bf Proposition}}[section]
\newtheorem{cor}{{\bf Corollary}}[section]

\theoremstyle{remark}

\theoremstyle{definition}

\author{Jamel BENAMEUR}
\address{Facult\'e des Sciences de Bizerte, D\'epartement
de Math\'ematiques\\ Zarzouna 7021, Tunisia.} \email{\it
jamel.benameur@fsb.rnu.tn}


\title{ASYMPTOTIC ANALYSIS OF MHD SYSTEMS}
\date{\today}

\begin{document}
\begin{abstract}In this paper, we
study the convergence of strong solutions of a Magneto-Hydro-Dynamic
system. On the torus ${\mathbb{T}}^3$, the proof is based on
Schochet's methods, whereas in the case of the whole space ${\bf
\mathbb{R}^3}$, we use Strichartz's type estimates and a product
law's $2D\times3D$.
\end{abstract}
\maketitle





\section{Introduction}
\noindent In this paper we study unique local existence and
asymptotic
 behaviour of solutions for the $\bf{MHD}$ system, see \cite{DDG},
$$\left\{\begin{array}{rcl} \displaystyle
\partial _tu+u.\nabla u
-\frac{E}{\varepsilon}\Delta u+\frac{e\times
u}{\varepsilon}-\frac{\Lambda}{\varepsilon}\mbox{curl}(b)\times
e'-\frac{\Lambda\theta}{\varepsilon}\mbox{curl}~b\times
b&=&-\displaystyle\frac{\nabla p}{\varepsilon}\\
\partial _tb-\displaystyle\frac{\Delta
b}{\theta}+u. \nabla b-b. \nabla u-
\mbox{curl}~\displaystyle\frac{u\times
e'}{\theta}&=&0\\
\mbox{div}~u=0,\;\mbox{div}~b&=&0
\end{array}\right.$$
in $\R^+_t\times\Omega_x$, where $\Omega$ is the whole space
$\mathbb{R}^3$ or the torus
 ${\mathbb{T}}^3:=\displaystyle{\R^3}/{\Z^3}$; $u$
is the velocity field, $b$ is the magnetic field and $e$, $e'$ are
two fixed vectors. The parameters $E$, $\varepsilon$, $\Lambda$ and
$\theta$ represent consecutively the Ekman number, the Rossby
number, the Elsasser number and the magnetic Reynolds number. We
notice that these parameters satisfy, according to \cite{DDG},
$$\varepsilon\rightarrow 0,\quad\Lambda=\mathcal{O}(1),
\quad\varepsilon\theta\rightarrow0\quad\quad\hbox{and}\quad\quad
E\sim\varepsilon^2.$$ These equations modalize the
magneto-hydro-dynamic flow in the Earth's core which is believed to
support a self-excited dynamo process generating the Earth's
magnetic field.  Here, we present the analytical study of a
simplified problem  where we choose $E=\varepsilon^2$,
$\Lambda=\varepsilon^{3/2}$, $\theta=\varepsilon^{-1/2}$ and
$e=e'=-e_3$. Then, we can see that our system take the following
form
$$\left\{
\begin{array}{rcl}
\partial _tu -\varepsilon\Delta u+
u.\nabla u-\mbox{curl } (b)\times b+ \sqrt{\varepsilon}\mbox{curl
}(b)\times e_3+\displaystyle\frac{u\times
e_3}{\varepsilon}&=&-\nabla p\\
\partial_t b- \sqrt{\varepsilon}\Delta b+
u. \nabla b-b. \nabla u+\sqrt{\varepsilon}\mbox{curl
}(u\times e_3)&=&0\\
\mbox{div}~u=0,\;\mbox{div}~b&=&0\\
\end{array}\right.\leqno{({\bf MHD^\varepsilon})}$$
The goal of this study is to find the "limit system" when
$\varepsilon$ goes to zero.\\
We denote by $\mathbb{P}$ the $L^2$ orthogonal projection on
divergence-free vector fields. Applying $\mathbb{P}$ to the first
equation of ${\bf (MHD^\varepsilon)}$, then one can see that
$U^\varepsilon=(u^\varepsilon,b^\varepsilon)$ is a solution of the
following abstract form
$$ \left\{
\begin{array}{rcl}
\partial_tU+Q(U, U)+
a_2^{\varepsilon}(D)U+L^{\varepsilon}(U)&=&0\quad in\quad
\mathbb{R}_t^+\times\Omega_x\\
\mbox{div}~u=0,\;\mbox{div}~b&=&0,\\
\end{array}
\right. \leqno{({\mathcal S}^{\varepsilon})} $$ where the
quadratic term $Q$ is defined by\\\\
 $$ Q(U,
U)=\Big({\mathbb{P}}(u.\nabla u)-{\mathbb{P}}(b.\nabla b), u. \nabla
b-b. \nabla u\Big), $$
the viscous term is $$
a_2^{\varepsilon}(D)U=\Big(-\varepsilon\Delta u,
-\sqrt{\varepsilon}\Delta b\Big),$$ and the linear perturbation
$L^{\varepsilon}$ is given by
 $$L^{\varepsilon}(U):=\Big({\mathbb{P}}(\frac{u\times
 e_3}{\varepsilon})+\sqrt{\varepsilon}\partial_3b, \sqrt{\varepsilon}\partial_3u\Big).$$
 In the literature, $a_2^{\varepsilon}(D)$ is elliptic and, in the case of a rotating
fluid $L^{\varepsilon}=\frac{1}{\varepsilon}L$ with $L$ a
skew-symmetric linear operator.\\
 Singular limits in system such as $({\mathcal S}^{\varepsilon})$ have been studied by
several authors. In the hyperbolic case, namely
$a_2^{\varepsilon}(D)U=0$, A. Babin, A. Mahalov and B. Nicolaenko
\cite{BMN} studied the incompressible rotating Euler equation on the
torus. Using the method introduced by S. Schochet (see \cite{S1} and
\cite{S2}), I. Gallagher studied in \cite{G1} this problem in its
abstract hyperbolic form. In the case of the incompressible rotating
Navier-Stokes equation on the torus, it is shown (see \cite{BMN} and
\cite{Gr1}) that the solutions converge to a solution of a certain
diffusion equation. Moreover, for a special initial condition, there
exists a sequence of solutions convergent to a solution of a two
dimensional Navier-Stokes equation. Motivated by this case, J. Y.
Chemin, B. Desjardin, I. Gallagher and E. Grenier studied in
\cite{CDGG} the incompressible fluids with anisotropic viscosity on
the whole space, the key of their proof is an anisotropic version of
Strichartz estimates. We refer to I. Gallagher \cite{G} for the
study of the abstract parabolic form. Among others, we also refer to
the basic results of \cite{CDGG}, \cite{G}, \cite{Gr1}, \cite{JMR1},
\cite{JMR2}, \cite{KM}, \cite{M}.\\Notice that the existence results
follow directly from the Friedrichs's method and the energy
estimates. Following the approximation scheme of Friedrichs, we
shall prove in section two global existence of ``Leray's solutions"
and local existence of strong solutions on uniform time; namely
solutions defined by the following results.
\begin{th1}\label{t41}
Let $\varepsilon> 0$ and $U_0=(u_0,b_0)\in L^2(\Omega)$ such that
$$\mbox{div}~u_0=0,\;\;\mbox{div}~b_0=0.$$ There exists
$U^\varepsilon:=(u^\varepsilon,b^\varepsilon)$ a solution of ${\bf
(MHD^\varepsilon)}$ with $U^\varepsilon\in L^\infty
(\mathbb{R}^+,L^2)\cap L^2 (\mathbb{R}^+,\dot{H}^1)$. Moreover,
$U^\varepsilon$ satisfies the following energy estimate. For all
$t\geq 0$,
\begin{eqnarray}
\|U^{\varepsilon}(t)\|_{L^2(\Omega)}^2+2\varepsilon\int_0^t\|\nabla
u^{\varepsilon}\|_{L^2}^2+2\sqrt{\varepsilon}\int_0^t\|\nabla
b^{\varepsilon}\|_{L^2}^2\leq\|U_0\|_{L^2}^2.
\end{eqnarray}
\end{th1}
\noindent \begin{th1}\label{t42} Let $U_0=(u_0,b_0)\in H^s(\Omega)$
with $s>\frac{3}{2}+2$ an integer, such that
$$\mbox{div}~u_0=0,\;\;\mbox{div}~b_0=0.$$ Then there exists $T> 0$, and a
constant $C> 0$ such that for all
  $\varepsilon\in ]0,1[ $, there exists a unique solution
$U^\varepsilon:=(u^\varepsilon,b^\varepsilon)\in{\mathcal
C}_T^0(H^s)\cap L^2_T(H^{s+1}) $ of system ${\bf (MHD^\varepsilon)}$
satisfying~; for all $t\in[0,T]$
\begin{equation}\label{eef}
\|U^{\varepsilon}(t)\|_{H^s(\Omega)}^2+2\varepsilon\int_0^t\|\nabla
u^{\varepsilon}\|_{H^s(\Omega)}^2+2\sqrt{\varepsilon}\int_0^t\|\nabla
b^{\varepsilon}\|_{H^s(\Omega)}^2\leq 2\|U_0\|_{H^s(\Omega)}^2.
\end{equation} Moreover, if $\|U_0\|_{H^s}\leq c\varepsilon$
$(c:=\frac{1}{C})$, then the solution is global.
\end{th1}
\noindent We are interested now to seeing the limit of strong
solutions of the $({\bf MHD})$ system when $\varepsilon$ goes to
zero. Since $\partial_t u^{\varepsilon}$ is not bounded in
$\varepsilon$, one cannot take the limit directly in the system, and
then the classical
proofs (see for example \cite{L1},\cite{T}) no longer work.\\\\
$\bullet$ In the case ${\bf \Omega=\mathbb{T}^3}$, a method to get
round this difficulty is to use the group ${\mathcal L}(t)$
associated with the operator $L(u):=\mathbb{P}(u\times e_3)$. We
consider the filtered solution $v^{\varepsilon}:={\mathcal
L}(-\frac{t}{\varepsilon})u^{\varepsilon}$ and we look for the limit
system in ${\mathcal D}'$ satisfied by the "eventuel" limit of
$(v^{\varepsilon}, b^{\varepsilon})$, $$\leqno({\bf
LS})\quad\quad\quad\quad\quad\quad\left\{
\begin{array}{clll}
\partial _tv+q^0(v,v)-\overline{{\mathbb{P}}(b.\nabla b)}&=&0
\\
\partial_tb+{\bar v}.\nabla b-b.\nabla {\bar v}&=&0\quad in\quad \mathbb{R}_t^+\times {\mathbb
T}_x^3\\ (v,b)_{/t=0}&=&(u_0,b_0)\\
\mbox{div}~v=\mbox{div}~b&=&0,\\
\end{array}
\right. $$ where $$q^0(v,v):=\lim_{\varepsilon}{\mathcal
L}(-\frac{t}{\varepsilon})\Big[{\mathbb{P}\Big(({\mathcal
L}(-\frac{t}{\varepsilon})v).\nabla ({\mathcal
L}(-\frac{t}{\varepsilon})v)\Big)}\Big]\quad in\quad {\mathcal D}'
$$ and here, we have denoted $\bar{f}(x_1,x_2)=\displaystyle{\int}
f(x_1,x_2,x_3)dx_3$ and $f_{osc}=f-\bar{f}$.\\ The system $({\bf
LS})$ is taking in the sobolev space, precisely we have the
following result.
\begin{th1}\label{t43}
Let $s>\frac{3}{2}+2$ be an integer and $U_0=(u_0,b_0)\in
H^s(\mathbb{T}^3)$ such that $\mbox{div}~u_0=0,\;\mbox{div}~b_0=0$.
Then there exists $T>0$, and a unique solution $(v,b)\in{\mathcal
C}_T^0(H^s)$ of system ${\bf (LS)}$.
\end{th1}
\noindent More, precisely we have the following convergence result.
\begin{th1}\label{t44}
Let $s>\frac{3}{2}+2$ be an integer and $U_0=(u_0,b_0)\in
H^s(\mathbb{T}^3)$ such that $\mbox{div}~u_0=0,\;\mbox{div}~b_0=0$.
We denote by $U^{\varepsilon}=(u^{\varepsilon},b^{\varepsilon})$ the
family of solutions of $({\bf MHD})$ given by Theorem \ref{t42}.
Then, for all $s'< s$,
$$\begin{array}{cllll}u^\varepsilon-{\mathcal
L}(\frac{t}{\varepsilon})v&=&o(1)\quad in\quad
 L^\infty_T(H^{s'}(\mathbb{T}^3))\\
 b^{\varepsilon}-b&=&o(1)\quad in\quad
 L^\infty_T(H^{s'}(\mathbb{T}^3)),\\
\end{array}$$
where $(v,b)$ is the solution of $({\bf LS})$. \end{th1} \noindent
$\bullet$ In the case ${\bf \Omega=\mathbb{R}^3}$, we have the
following result
\begin{th1}\label{t45}
Let $s>\frac{3}{2}+2$ be an integer and $U_0=(u_0,b_0)\in
H^s(\mathbb{R}^3)$ such that $\mbox{div}~u_0=0,\;\mbox{div}~b_0=0$.
We denote by $U^{\varepsilon}=(u^{\varepsilon},b^{\varepsilon})$ the
family of solutions of $({\bf MHD})$ given by Theorem \ref{t42}.
Then, for all $s'< s$
$$\begin{array}{cllll}u^\varepsilon&=&o(1)\quad in\quad
 L^4_T({\mathcal C}^{s'-\frac{3}{2}}(\mathbb{R}^3))\\
 b^{\varepsilon}-b_0&=&o(1)\quad in\quad
 L^{\infty}_T(H^{s'}(\mathbb{R}^3)).\\
\end{array}$$
\end{th1}
\noindent Since the limit of the system $({\bf MHD})$ is the
bidimensional Navier-Stokes equations, it is natural to consider
initial data of the type $(u_0,b_0)=({\bar u}_0+w_0,b_0)$, where
${\bar u}_0={\bar u}_0(x_h)$ and $w_0=w_0(x_h,x_3),
b_0=~b_0(x_h,x_3)$ (see \cite{CDGG}) .\\ Before stating the results,
it will be useful to consider the following system
$$ \left\{
\begin{array}{cllll}
\partial_t{\bar u}-\varepsilon\Delta{\bar u}+
{\bar u}.\nabla_h\;{\bar u}&=&-(\nabla_h
p,0)\quad in\quad \mathbb{R}_t^+\times\mathbb{R}_h^2\\
\mbox{div}_h~{\bar u}&=&0\\ {\bar u}_{/t=0}&=&{\bar u}_0.
\end{array}
\right. \leqno({\bf NS2D^\varepsilon})$$ Using the classical
Friedrich's scheme, we can prove the existence of strong solutions
on uniform time for the system $({\bf NS2D^\varepsilon})$.
Precisely, we have the classical result
\begin{th1}\label{t46}
Let ${\bar u}_0 \in \Big(H^{\sigma}(\mathbb{R}^2)\Big)^3$ be a
divergence-free vector field with $\sigma>\frac{2}{2}+2$ an integer.
Then there exists $T_0:=\displaystyle\frac{1}{C(\sigma){\|\bar
u}_0\|_{H^\sigma}}$ such that for all
  $\varepsilon> 0 $, there exists a unique solution
${\bar u}^\varepsilon\in{\mathcal C}^0_{T_0}
H^\sigma(\mathbb{R}^2))\cap L^2_{T_0}( H^{\sigma+1}(\mathbb{R}^2)) $
of system $({\bf NS2D^\varepsilon})$ satisfying; for all
$t\in[0,T_0]$ $$ \|{\bar
u}^{\varepsilon}(t)\|_{H^\sigma}^2+2\varepsilon\int_0^t\|\nabla
{\bar u}^{\varepsilon}\|_{H^\sigma}^2\leq 2\|{\bar
u}_0\|_{H^\sigma}^2. $$
\end{th1}
\noindent We pose $$\left\{
\begin{array}{cllll}
 (u_0,b_0)&=&({\bar
u}_0+w_0,b_0)\\ \mbox{div}_h~{\bar
u}_0&=&\mbox{div}~w_0=\mbox{div}~b_0=0\\ {\bar u}_0&=&{\bar
u}_0(x_h)\in H^{s+1}(\mathbb{R}^2)\\ w_0&=&w_0(x_h,x_3)\in
H^{s}(\mathbb{R}^3)\\ b_0&=&b_0(x_h,x_3)\in H^{s}(\mathbb{R}^3) \\
\mbox{with}\quad&& s>\frac{3}{2}+2 \quad\mbox{an integer}.
\end{array} \right. \leqno({\bf 2D\times3D})$$ We suppose the
condition is satisfy and we denoted $({\bar u}^\varepsilon)$ the
family of solutions of $({\bf NS2D})$.\\ Now we considered the
following system
 $$ \left\{
\begin{array}{cllll}
\partial _tw-\varepsilon\Delta w
+\displaystyle\frac{1}{\varepsilon}w\times e_3
+\sqrt{\varepsilon}\partial_3B-B.\nabla
B&=&-\nabla p_L\\
\partial_tB- \displaystyle{\sqrt{\varepsilon}}\Delta B
+\sqrt{\varepsilon}\partial_3w+({\bar u}^\varepsilon+w).\nabla
B&&\quad in\quad\quad\quad
[0,T_0]\times\mathbb{R}^3\\
\quad\quad\quad\quad\quad\quad\quad\quad\quad\quad\quad
-B.\nabla ({\bar u }^\varepsilon+w)&=&0\\
(w,B)(0)&=&(w_0,b_0)\\
(\mbox{div}~w,\mbox{div}~B)&=&(0,0).\\
\end{array}
\right. \leqno({\bf MHD3D^\varepsilon})$$
\begin{th1}\label{t47}Suppose that the condition $({\bf
2D\times3D})$ is satisfied. There exists $0<T_1\leq T_0$ such that
for all
  $\varepsilon> 0 $, there exists a unique solution
$(w^\varepsilon,B^\varepsilon)\in L^\infty_{T_1}(
H^s(\mathbb{R}^3))\cap L^2_{T_1}(H^{s+1}(\mathbb{R}^3)) $ of system
${\bf (MHD3D^\varepsilon)}$ satisfying ; for all $t\in[0,T_1]$ $$
\|w^{\varepsilon}(t)\|_{H^s}^2+\|B^{\varepsilon}(t)\|_{H^s}^2+
2\varepsilon\int_0^t\|\nabla
w^{\varepsilon}\|_{H^s}^2+2\sqrt{\varepsilon}\int_0^t\|\nabla
B^{\varepsilon}\|_{H^s}^2\leq 2\|(w_0,b_0)\|_{H^s}^2. $$ Morowever
$$w^\varepsilon\rightarrow 0\quad\mbox{in}\quad
L^4([0,T_1],{\mathcal C}^{s'-\frac{3}{2}}({\mathbb R}^3))~;~\forall
s'< s.$$
\end{th1}
\noindent
\begin{th1}\label{t48}
Suppose that the condition $({\bf 2D\times3D})$ is satisfy.\\Then,
for all $\varepsilon> 0$, there exists a unique solution
$(u^\varepsilon,b^\varepsilon)$ of the system $({\bf
MHD^\varepsilon})$ such that $$ u^\varepsilon-{\bar
u}^\varepsilon,\quad b^\varepsilon \in
L^\infty_{T_1}(H^s(\mathbb{R}^3))\cap
L^2_{T_1}(H^{s+1}(\mathbb{R}^3)).$$ Moreover, we have for all $t\in
[0,T_1]$,
$$ \|u^\varepsilon(t)-{\bar
u}^\varepsilon(t)\|_{H^s}^2+\|b^\varepsilon(t)\|_{H^s}^2+
2\varepsilon\displaystyle\int_0^t\|\nabla (u^\varepsilon-{\bar
u}^\varepsilon)(\tau)\|_{H^s}^2\;d\tau+2\sqrt{\varepsilon}\displaystyle\int_0^t\|\nabla
b^{\varepsilon}\|_{H^s}^2\;d\tau\leq 2\|(w_0,b_0)\|_{H^s}^2.$$
\end{th1}
\noindent Now we are ready to state the main convergence result in
the case of the whole space $\mathbb{R}^3$.
\begin{th1}\label{t49}
We keep the same hypothesis as in Theorem \ref{t48} above and we
suppose $s>\frac{3}{2}+4$ an integer. Then, for all $s'> s$
$$ u^\varepsilon-{\bar u}^\varepsilon-w^\varepsilon =
o(1)\quad\mbox{in}\quad L^\infty_{T_1}(H^{s'}(\mathbb{R}^3))$$
$$b^\varepsilon-B^\varepsilon= o(1)\quad\mbox{in}\quad
L^\infty_{T_1}(H^{s'}(\mathbb{R}^3)).$$
\end{th1}
\noindent The structure of this paper is as follows. In the next
section, we present the proofs of the existence theorems (Theorem
\ref{t41}, \ref{t42}). The third section is devoted to the proof of
the convergence result in the case $\Omega=\mathbb{T}^3$ (Theorem
\ref{t44}) and the study of the system $(\bf LS)$ (Theorem
\ref{t43}). In the final section, we consider the case of the whole
space $\mathbb{R}^3$; We give the proof of the Theorems \ref{t45},
\ref{t47}, \ref{t48}, \ref{t49}.
\section{Existence results}
\subsection{Proofs of Theorems \ref{t41} and \ref{t42}}
In this section we shall prove Theorems \ref{t41} and \ref{t42}.
\\ We begin by
observing that using the energy methods, one can prove global
existence of so-called ``Leray's solutions" for the system $({\bf
MHD}^\varepsilon)$. The crucial fact is the following $L^2$-energy
estimate
\begin{equation}
\label{Leray}
\|U^\varepsilon(t)\|_{L^2}^2+2_\varepsilon\int_0^t\|\nabla
u^\varepsilon(\tau)
\|_{L^2}^2\;d\tau+2\sqrt{\varepsilon}\int_0^t\|\nabla
b^\varepsilon(\tau) \|_{L^{2}}^2\;d\tau\leq \|U_0\|_{L^2}^2.
\end{equation}
We now turn to the case of strong solutions. Let us introduce, for a
strictly positive integer $n$, the Friedrich's operator $J_n$
defined by: $$ J_n\;u={\mathcal F}^{-1}\Big({\bf
1}_{B(0,n)}{\mathcal F}(u)\Big).$$ After this definition, we
consider the following approximate Magneto-Hydro-Dynamic system
$({\bf MHD}_n)$ $$ \left\{
\begin{array}{cllll}
\partial _tu_n -\varepsilon \Delta J_n u_n+
J_n\mbox{div}~(J_n u_n\otimes J_n u_n)-J_n\mbox{div}~(J_n b_n\otimes
J_n b_n)+\sqrt{\varepsilon}\partial_3 (J_n b_n)+\\
\displaystyle\frac{ J_n u_n\times e}{\varepsilon} =
\nabla{\Delta}^{-1}\mbox{div}~\Big( J_n\mbox{div}~(J_n u_n\otimes
J_n u_n)-J_n\mbox{div}~(J_n b_n\otimes J_n b_n)+\frac{ J_n
u_n\times e}{\varepsilon}\Big),\\
\partial _t b_n - \sqrt{\varepsilon}\Delta J_n b_n+J_n\mbox{div}~(J_n u_n\otimes J_n b_n)-J_n\mbox{div}~(J_n b_n\otimes
J_n u_n)+\sqrt{\varepsilon}\partial_3 (J_n u_n)=0,\\ (u_n|_{t=0},
b_n|_{t=0})=(J_n u_0, J_n b_0).
\end{array}
\right.$$ By the theory of ordinary differential equations in $H^s$
we know that the system $({\bf MHD}_n)$ has a unique maximal
solution $U_n:=(u_n,b_n)$ in the space ${\mathcal C}^1([0,
T_n^*(\varepsilon)[, H^s)$. Using uniqueness and the fact that
$\mbox{div}~u_n=\mbox{div}~b_n=0$ and $J_n^2=J_n$ we can re-write
the system $$ \left\{
\begin{array}{cllll}
\partial _tu_n -\varepsilon \Delta u_n+
J_n( u_n.\nabla u_n)-J_n( b_n.\nabla
b_n)+\sqrt{\varepsilon}\partial_3 b_n+\displaystyle\frac{
 u_n\times e}{\varepsilon} = \\ \nabla{\Delta}^{-1}\mbox{div}~\Big(
J_n\mbox{div}~( u_n\otimes u_n)-J_n\mbox{div}~( b_n\otimes
b_n)+\displaystyle\frac{ u_n\times e}{\varepsilon}\Big),\\
\partial _t b_n - \sqrt{\varepsilon}\Delta b_n+J_n(u_n.\nabla b_n)-
J_n(b_n.\nabla u_n)+\sqrt{\varepsilon}\partial_3 u_n=0,\\
(u_n|_{t=0}, b_n|_{t=0})=(J_n u_0, J_n b_0).
\end{array}
\right.\leqno({\bf MHD}_n)$$ To continue the proof, we recall
without proof the following product law.
\begin{lem}\label{l41}
Let $\sigma>\frac{3}{2}+2$ be an integer and $a$, $b$ and $c$ three
vectors field in $H^{\sigma}(\Omega)$ such that $\mbox{div}~a=0$.
Then, a constant $C$ exists such that $$
\begin{array}{cllll}
|{< a.\nabla b,b>}_{H^\sigma}|&\leq  & C\|\nabla a
\|_{H^{\sigma-1}}\|\nabla b \|_{H^{\sigma-1}}^2\\ |{< a.\nabla
b,c>}_{H^\sigma}+{< a.\nabla c,b>}_{H^\sigma}|&\leq & C\|\nabla
a\|_{H^{\sigma-1}}\|\nabla b \|_{H^{\sigma-1}}\|\nabla
c\|_{H^{\sigma-1}}.
\end{array}
$$
\end{lem}
\noindent (See the Appendix for the proof of the lemma.)\\\\
 We take the scalar
product in $H^s$ and we use the lemma above, we obtain for all
$t\in[0,T_n^*(\varepsilon)[$,
\begin{equation}\label{Fr0}
\frac{1}{2}\frac{d}{dt}\|U_n(t)\|_{H^s}^2+\varepsilon\|\nabla u_n(t)
\|_{H^{s}}^2+\sqrt{\varepsilon}\|\nabla b_n(t) \|_{H^{s}}^2\leq
 C\|\nabla U_n(t) \|_{H^{s-1}}^3.
\end{equation}
Then
\begin{equation}\label{Fr1}
\|U_n(t)\|_{H^s}^2+2\varepsilon\int_0^t\|\nabla u_n(\tau)
\|_{H^{s}}^2d\tau+2\sqrt{\varepsilon}\int_0^t\|\nabla b_n(\tau)
\|_{H^{s}}^2d\tau\leq
 \|U_0\|_{H^s}^2+C \int_0^t\|\nabla U_n(\tau) \|_{H^{s-1}}^3d\tau.
\end{equation}
We set $T(n,\varepsilon):=\mbox{Sup}\{0\leq t<
T_n^*(\varepsilon);\quad \forall
\tau\in[0,t],\|U_n(\tau)\|_{H^s}\leq 2\|U_0\|_{H^s}\}$. Using
(\ref{Fr1}) and Gronwall lemma we obtain, for all
$t\in[0,T(n,\varepsilon)[$ , $$\|U_n(t)\|_{H^s}^2\leq
\|U_0\|_{H^s}^2\exp(2Ct\|U_0\|_{H^s}).$$ Thus,
$$T(n,\varepsilon)> T:=\frac{1}{C\|U_0\|_{H^s}}> 0.$$
Moreover, for all $t\in[0,T]$,
\begin{equation}  \quad\quad \quad
\|U_n(t)\|_{H^s}^2+2\varepsilon\int_0^t\|\nabla u_n(\tau)
\|_{H^{s}}^2\;d\tau+2\sqrt{\varepsilon}\int_0^t\|\nabla b_n(\tau)
\|_{H^{s}}^2\;d\tau\leq 2\|U_0\|_{H^s}^2.\label{Fr2}
\end{equation}
Now, the problem is to pass to the limit. Using Ascoli's theorem,
the Cantor's diagonal process as in Navier-Stokes equations (see
\cite{Ch}) and the estimate (\ref{Fr2}), we obtain a solution
satisfying, for all $t\in[0,T]$,
\begin{equation}
\label{Fr3}
\|U^\varepsilon(t)\|_{H^s}^2+2\varepsilon\int_0^t\|\nabla
u^\varepsilon(\tau)
\|_{H^{s}}^2\;d\tau+2\sqrt\varepsilon\int_0^t\|\nabla
b^\varepsilon(\tau) \|_{H^{s}}^2\;d\tau\leq 2\|U_0\|_{H^s}^2.
\end{equation}
This regularity implies in a standard way the uniqueness. It remains
to prove the global existence when the initial data is small enough.
We assume now that $\|U_0\|_{H^s}\leq c\varepsilon$,\quad
($c=\frac{1}{C}$), and we set $$ T_n(\varepsilon):=\mbox{Sup}\{0\leq
t< T_n^*(\varepsilon);\quad \forall
\tau\in[0,t],\|U_n(\tau)\|_{H^s}\leq c\varepsilon\}.$$ Using
(\ref{Fr1}), it suffices to show that $T_n(\varepsilon)=
T_n^*(\varepsilon)$. By (\ref{Fr0}) we have
$\displaystyle\frac{d}{dt}\|U_n\|_{H^s}^2(0)< 0$, then there exists
$t_n> 0$ such that $\|U_n(t_n)\|_{H^s}< c\varepsilon.$ Since the
quantity $\|U_n(t)\|_{H^s}$ is decreasing on
$[t_n,T_n^*(\varepsilon)[$, then
$T_n(\varepsilon)=T_n^*(\varepsilon)$. This achieves the
proof.\endproof
\section{The case $\Omega=\mathbb{T}^3$}
\noindent Let $(U^{\varepsilon})$ be a family of strong solutions of
the system (${\mathcal S}^\varepsilon$) with initial data $U_0$. To
take the limit when $\varepsilon\longrightarrow 0$, the classical
proofs no longer work because $(\partial_t u^{\varepsilon})$ is not
uniformly bounded. An idea (as in \cite{S2} for instance) is to
``filter" the system by the group ${\mathcal L}(t)$ associated to
$L$.\\ In what follows, we recall some properties of the Coriolis
force $L(u)$. We consider, as in \cite{Gr1}, the "wave equation" $$
\left\{
\begin{array}{cllll}
&\partial_tu+L(u)=0\quad in\quad \R_t\times{\mathbb{T}}_x^3,\\
&u(0)=u_0\quad\mbox{with}\quad \mbox{div}~u_0=0.\\
\end{array}
\right. $$
\begin{lem}\label{l42}
The above system has a global solution denoted by $u(t)={\mathcal
L}(t)u_0$, such that for all  $s\in \R$ and for all $u_0\in
H^s({\mathbb{T}}^3)$, $$ \|{\mathcal
L}(t)u_0\|_{H^s({\mathbb{T}}^3)}=\|u_0\|_{H^s({\mathbb{T}}^3)}\quad
and\quad \|^t{\mathcal
L}(t)u_0\|_{H^s({\mathbb{T}}^3)}=\|u_0\|_{H^s({\mathbb{T}}^3)} .$$
\\ Moreover, if we denote by $k=(k_1,k_2,k_3)$ the Fourier coordinates,
then $u$ is explicitly given by $$ {\mathcal
F}u(t,k)=\exp(i\omega(k)t)({\mathcal F}u(0,k),\nu_k^+)\nu_k^+
+\exp(-i\omega(k)t)({\mathcal F}u(0,k),\nu_k^-)\nu_k^- ,$$ where
$\omega(k)=\frac{k_3}{| k|}$, $\nu_k^{\pm}$ are given unit vectors
and $(.,.)$ denotes the usual scalar product.
\end{lem}
Now we define $$ v^{\varepsilon}(t) = {\mathcal
L}\Big(\frac{-t}{\varepsilon}\Big)u^{\varepsilon}(t)$$ then,
$V^{\varepsilon}:=(v^{\varepsilon}, b^{\varepsilon})$, satisfies the
following system $$ \left\{
\begin{array}{cllll}
&\partial _tv^{\varepsilon} -\varepsilon \Delta
v^{\varepsilon}+{\mathcal
L}\Big(\frac{-t}{\varepsilon}\Big){\mathbb{P}}({\mathcal
L}\Big(\frac{t}{\varepsilon}\Big)v^{\varepsilon}.\nabla {\mathcal
L}\Big(\frac{t}{\varepsilon}\Big)v^{\varepsilon})
{\mathcal L}\Big(\frac{-t}{\varepsilon}\Big){\mathbb{P}}(\mbox{curl}
b^{\varepsilon}\times b^{\varepsilon}) +\sqrt{\varepsilon}{\mathcal
L}\Big(\frac{-t}{\varepsilon}\Big){\mathbb{P}}(\mbox{curl}
b^{\varepsilon}\times e_3)=0\\ &\partial_tb^{\varepsilon}
-\sqrt{\varepsilon} \Delta b^{\varepsilon}+({\mathcal
L}\Big(\frac{t}{\varepsilon}\Big)v^{\varepsilon}). \nabla
b^{\varepsilon}-b^{\varepsilon}. \nabla ({\mathcal
L}\Big(\frac{t}{\varepsilon}\Big)v^{\varepsilon})
+\sqrt{\varepsilon}\mbox{curl}(({\mathcal
L}\Big(\frac{t}{\varepsilon}\Big)v^{\varepsilon})\times e_3)=0\\
&\quad\quad\quad\quad\quad\quad\quad\quad
\mbox{div}~v^{\varepsilon}=\mbox{div}~b^{\varepsilon}=0\\
&\quad\quad\quad\quad\quad\quad\quad\quad (v^{\varepsilon},
b^{\varepsilon})(0)=(u_0, b_0).\\
\end{array}
\right. $$ This system can be re-written in the following way $$
\left\{
\begin{array}{cllll}
&\partial_tV^{\varepsilon}+Q^{\varepsilon}(V^{\varepsilon},
V^{\varepsilon})+ a_2^{\varepsilon}(D)V^{\varepsilon}+{\tilde
L}^{\varepsilon}(V^{\varepsilon})=0\quad in\quad
\R^+_t\times{\mathbb{T}}_x^3\\ &\mbox{div}~V^{\varepsilon}=0\\
&V^{\varepsilon}(0)=U_0=(u_0, b_0),\\
\end{array}
\right. \leqno{({\tilde{\mathcal S}}^{\varepsilon})} $$ where the
``filtered" quadratic form $Q^\varepsilon$ is given by
$$
Q^{\varepsilon}(V, V)=\Big({\mathcal
L}\Big(\frac{-t}{\varepsilon}\Big){\mathbb{P}}({\mathcal
L}\Big(\frac{t}{\varepsilon}\Big)v.\nabla {\mathcal
L}\Big(\frac{t}{\varepsilon}\Big)v)-{\mathcal
L}\Big(\frac{-t}{\varepsilon}\Big){\mathbb{P}}(b.\nabla
b),({\mathcal L}\Big(\frac{t}{\varepsilon}\Big)v). \nabla b-b.
\nabla ({\mathcal L}\Big(\frac{t}{\varepsilon}\Big)v)\Big),
$$and, $$ a_2^{\varepsilon}(D)V=\Big(-\varepsilon \Delta v,
-\sqrt{\varepsilon}\Delta b\Big), $$ $$ {\tilde
L}^{\varepsilon}(V)=\sqrt{\varepsilon}\Big({\mathcal
L}\Big(\frac{-t}{\varepsilon}\Big)\partial_3 b,
\partial_3{\mathcal L}\Big(\frac{t}{\varepsilon}\Big)v\Big) . $$
When $\varepsilon$ goes to 0, we obtain formally the following limit
system $$
 \left\{
\begin{array}{cllll}
&\partial_t V+Q^0(V,V)=0\quad in\quad
\R^+_t\times{\mathbb{T}}_x^3\\ &\mbox{div}~v=\mbox{div}~b=0\\
&V(0)=U_0=(u_0, b_0),\\
\end{array}
\right.\leqno{({\bf LS})} $$ where $Q^0(V,V)$ is the limit in
${\mathcal D}^{'}$ of $Q^{\varepsilon}(V,V)$.
\subsection{Proof of Theorem \ref{t43}}
The proof is similar to the one of Theorem \ref{t42}. We have just
to estimate the term $$\Big|\int_0^t
<q^0(v_n,v_n),v_n>_{H^s}\Big|\;,
$$ where $(v_n,b_n)$ is the solution of the approximate limit
system. Observe that $$\lim_{\varepsilon}\int_0^t
<q^\varepsilon(v_n,v_n),v_n>_{H^s}=\int_0^t
<q^0(v_n,v_n),v_n>_{H^s}$$ and using the product law given by Lemma
\ref{l41}, we obtain $$\Big|\int_0^t
<q^0(v_n,v_n),v_n>_{H^s}\Big|\leq C \int_0^t\|\nabla v_n(\tau)
\|_{H^{s-1}}^3\;d\tau$$ which completes the proof.\endproof
\subsection{Proof of Theorem \ref{t44}}
The proof of this theorem is based on a method used in \cite{G},
\cite{Gr1} for instance.

Let $W^{\varepsilon}= V^{\varepsilon}-V=(v^{\varepsilon}-v,
b^{\varepsilon}-b) = (W_1^{\varepsilon},W_2^{\varepsilon}),$ then
$W^{\varepsilon}$ satisfies
 $$
 \left\{\begin{array}{cllll}
&\partial_tW^{\varepsilon}+Q^{\varepsilon}(W^{\varepsilon},
W^{\varepsilon}+2V)+ a_2^{\varepsilon}(D)W^{\varepsilon}+{\tilde
L}^{\varepsilon}(W^{\varepsilon})= F^{\varepsilon} +
R_{osc}^{\varepsilon}\\
&\mbox{div}~W^{\varepsilon}_1=\mbox{div}~W^{\varepsilon}_2=0\\
&W^{\varepsilon}(0)=(0, 0),\\
\end{array}
\right.$$ where
\begin{eqnarray*}F^{\varepsilon}&=&-a_2^{\varepsilon}(D)V-{\tilde
L}^{\varepsilon}V\\ &=&(\varepsilon\Delta
v-\sqrt{\varepsilon}{\mathcal
L}(-\frac{t}{\varepsilon})b;\sqrt{\varepsilon}\Delta
b-\sqrt{\varepsilon}{\mathcal L}(\frac{t}{\varepsilon})v)\\
\end{eqnarray*}
 and
$R_{osc}^{\varepsilon}=Q^0(V, V)-Q^{\varepsilon}(V, V)$.\\ We recall
that $v$ and $b$ are in the space ${\mathcal C}_T^0( H^s)$, then
$$F^{\varepsilon}\longrightarrow 0\quad\mbox{ in}\quad {\mathcal
C}_T^0(H^{s-2})$$ precisely we are
\begin{equation}
\|F^{\varepsilon}\|_{{\mathcal C}_T^0( H^{s-2})}\leq
C\sqrt\varepsilon.
\end{equation}
 On the other hand, the oscillating term $R_{osc}^{\varepsilon}$
can be written as follows.
\begin{eqnarray}
\nonumber
 R_{osc}^{\varepsilon}&=&\Big({\mathcal
L}(-\frac{t}{\varepsilon}){\mathbb{P}}({\mathcal
L}(\frac{t}{\varepsilon})v.\nabla {\mathcal
L}(\frac{t}{\varepsilon})v)-q^0(v,v),0\Big)\nonumber \\
&&+\Big({\mathcal L}(-\frac{t}{\varepsilon}){\mathbb{P}}(b.\nabla
b)-{\overline{{\mathbb{P}}(b.\nabla b)}},0\Big)+\Big(0,-({\mathcal
L}(\frac{t}{\varepsilon})v).\nabla b+\overline{v}.\nabla b\Big)
\nonumber\\&&+\Big(0,b.\nabla\overline{v}-b.\nabla({\mathcal
L}(\frac{t}{\varepsilon})v)\Big), \nonumber
\end{eqnarray} where
 $$ q^0(v,v)=\lim_ {\varepsilon}{\mathcal
L}(-\frac{t}{\varepsilon}){\mathbb{P}}({\mathcal
L}(\frac{t}{\varepsilon})v.\nabla {\mathcal
L}(\frac{t}{\varepsilon})v)\quad in\quad{\mathcal D}'.$$In the
sequel, for any three-vector field $X$, we shall note, $$
X^{\pm}(n)=(X,\nu^{\pm}(n))\nu^{\pm}(n).$$
 We can write $$
R_{osc}^{\varepsilon}=\sum_{k=0}^3A_k^{\varepsilon}. $$
 We begin by studying the term $A_0^{\varepsilon}$.
\begin{eqnarray}
\nonumber A_0^{\varepsilon}&=&{\mathcal
L}(-\frac{t}{\varepsilon}){\mathbb{P}}({\mathcal
L}(\frac{t}{\varepsilon})v.\nabla {\mathcal
L}(\frac{t}{\varepsilon})v)-q^0(v,v)\\ \nonumber
&=&q^{\varepsilon}(v,v)-q^0(v,v).
\end{eqnarray}
By the non stationary phase's theorem we have
\begin{eqnarray} {\mathcal F}(A_0^{\varepsilon})(t,n)=\sum_{\sigma \in \{\pm\}^3}
\sum_{k+m=n\atop \omega_\sigma(n,k,m)\neq 0}
e^{-i\frac{t}{\varepsilon} \omega_\sigma(n,k,m)}( {({\mathcal
F}v)}^{\sigma_1}(t,k).{({\mathcal F}\nabla
v)}^{\sigma_2}(t,m))^{\sigma_3}(n),\nonumber\end{eqnarray} where
$\sigma=(\sigma_1,\sigma_2,\sigma_3)$ and
$\omega_\sigma(n,k,m)=\sigma_1\frac{n_3}{|n|}-\sigma_2\frac{k_3}{|k|}-\sigma_3\frac{m_3}{|m|}.$\\
Hence we can write $$A_0^{\varepsilon}= {\mathcal
F}^{-1}\Big(\sum_{\sigma \in \{\pm\}^3}\sum_{k+m=n\atop
\omega_\sigma(n,k,m)\neq 0} e^{-i\frac{t}{\varepsilon}
\omega_\sigma(n,k,m)}r_{\sigma}(n,k,m)
f_{\sigma}(t,k)g_{\sigma}(t,m)\Big). $$ \\ Similarly, we have $$
A_1^{\varepsilon}={\mathcal F}^{-1}\Big({\bf 1}_{\{n_3\neq
0\}}(n)e^{-i\frac{t}{\varepsilon}\frac{n_3}{| n|}}{({\mathcal
F}(b.\nabla b))}^+(t,n)+{\bf 1}_{\{n_3\neq
0\}}(n)e^{i\frac{t}{\varepsilon}\frac{n_3}{| n|}}{({\mathcal
F}(b.\nabla b))}^-(t,n)\Big),$$
 $$ A_2^{\varepsilon}={\mathcal
F}^{-1}\Big(\sum_{k_3\neq 0}e^{i\frac{t}{\varepsilon}\frac{k_3}{|
k|}}{({\mathcal F}v)}^+(t,k).({\mathcal F}\nabla
b)(t,n-k)+\sum_{k_3\neq 0}e^{-i\frac{t}{\varepsilon}\frac{k_3}{|
k|}}{({\mathcal F}v)}^-(t,k).({\mathcal F}\nabla b)(t,n-k)\Big)$$
and finally,
 $$
A_3^{\varepsilon}={\mathcal F}^{-1}\Big(\sum_{k_3\neq
0}e^{i\frac{t}{\varepsilon}\frac{k_3}{| k|}}({\mathcal
F}b)(t,n-k).{({\mathcal F}\nabla v)}^+(t,k)+\sum_{k_3\neq
0}e^{-i\frac{t}{\varepsilon}\frac{k_3}{| k|}}({\mathcal
F}b)(t,n-k).{({\mathcal F}\nabla v)}^-(t,k)\Big).$$ \\Since
$R_{osc}^{\varepsilon}$ converges weakly to 0 (but not
``strongly''), we shall divide it (as in \cite{G} for instance) into
high and low frequencies terms. Precisely, for any integer $N>1$, we
define
 $$A_{0,N}^{\varepsilon}={\mathcal
F}^{-1}\Big({\bf 1}_{\{| n|\leq N\}}\sum_{\sigma \in
\{\pm\}^3}\sum_{{k+m=n \atop
  \omega_{\sigma}(n,k,m)\neq0}\atop
  | k|,| m|\leq N}
e^{-i\frac{t}{\varepsilon}\omega_{\sigma}(n,k,m)}r_{\sigma}(n,k,m)
f_{\sigma}(t,k)g_{\sigma}(t,m)\Big), $$

$$
A_{1,N}^{\varepsilon}={\mathcal F}^{-1}\Big({\bf 1}_{\{n_3\neq 0,|
n|\leq N \}}(n)e^{-i\frac{t}{\varepsilon}\frac{n_3}{|
n|}}{({\mathcal F}(b.\nabla b))}^+(t,n)+{\bf 1}_{\{n_3\neq 0,|
n|\leq N\}}(n)e^{i\frac{t}{\varepsilon}\frac{n_3}{| n|}}{({\mathcal
F}(b.\nabla b))}^-(t,n)\Big),
$$
$$A_{2,N}^{\varepsilon}={\mathcal
F}^{-1}\Big({\bf 1}_{\{| n|\leq N\}} \displaystyle\sum_{k_3\neq 0
\atop | n-k|,| k|\leq N}e^{i\frac{t}{\varepsilon}\frac{k_3}{|
k|}}{({\mathcal F}v)}^+(t,k).({\mathcal F}\nabla
b)(t,n-k)+e^{-i\frac{t}{\varepsilon}\frac{k_3}{| k|}}{({\mathcal
F}v)}^-(t,k).({\mathcal F}\nabla b)(t,n-k)\Big)
$$
and$$ A_{3,N}^{\varepsilon}={\mathcal F}^{-1}\Big({\bf 1}_{\{|
n|\leq N\}} \sum_{k_3\neq 0 \atop
  | n-k|,| k|\leq N}e^{i\frac{t}{\varepsilon}\frac{k_3}{| k|}}({\mathcal
F}b)(t,n-k).{({\mathcal F}\nabla
v)}^+(t,k)+e^{-i\frac{t}{\varepsilon}\frac{k_3}{| k|}}({\mathcal
F}b)(t,n-k).{({\mathcal F}\nabla v)}^-(t,k)\Big). $$  Now, the idea
is to absorb the low frequency terms. For that, we set $$ {\tilde
A}_{0,N}^{\varepsilon}={\mathcal F}^{-1}\Big({\bf 1}_{\{| n|\leq
N\}} \sum_{\sigma \in \{\pm\}^3}\sum_{{k+m=n \atop
  \omega_{\sigma}(n,k,m)\neq0}\atop |k|,| n|\leq N}
\frac{e^{-i\frac{t}{\varepsilon}\omega_{\sigma}(n,k,m)}}{i\omega_{\sigma}(n,k,m)}
r_{\sigma}(n,k,m) f_{\sigma}(t,k)g_{\sigma}(t,m)\Big), $$
$$ {\tilde A}_{1,N}^{\varepsilon}={\mathcal
F}^{-1}\Big({\bf 1}_{\{n_3\neq 0,| n|\leq N
\}}(n)\frac{e^{-i\frac{t}{\varepsilon}\frac{n_3}{|
n|}}}{i\frac{n_3}{| n|}}({\mathcal F}(b.\nabla b))^+(t,n)+{\bf
1}_{\{n_3\neq 0,| n|\leq
N\}}(n)\frac{e^{i\frac{t}{\varepsilon}\frac{n_3}{|
n|}}}{-i\frac{n_3}{| n|}}({\mathcal F}(b.\nabla b))^-(t,n)\Big),
$$ $$
 {\tilde A}_{2,N}^{\varepsilon}={\mathcal
F}^{-1}\Big({\bf 1}_{\{| n|\leq N\}} \sum_{k_3\neq 0 \atop
  | n-k|,| k|\leq N}
\frac{e^{i\frac{t}{\varepsilon}\frac{k_3}{| k|}}}{-i\frac{k_3}{|
k|}}({\mathcal F}v)^+(t,k).({\mathcal F}\nabla
b)(t,n-k)+\frac{e^{-i\frac{t}{\varepsilon}\frac{k_3}{|
k|}}}{i\frac{k_3}{| k|}}({\mathcal F}v)^-(t,k).({\mathcal F}\nabla
b)(t,n-k)\Big)
$$and $$ {\tilde A}_{3,N}^{\varepsilon}={\mathcal F}^{-1}\Big({\bf
1}_{\{| n|\leq N\}} \sum_{k_3\neq 0 \atop |n-k|,| k|\leq
N}\frac{e^{i\frac{t}{\varepsilon}\frac{k_3}{| k|}}}{-i\frac{k_3}{|
k|}}({\mathcal F}b)(t,n-k).{({\mathcal F}\nabla
v)}^+(t,k)+\frac{e^{-i\frac{t}{\varepsilon}\frac{k_3}{|
k|}}}{i\frac{k_3}{| k|}}({\mathcal F}b)(t,n-k).{({\mathcal F}\nabla
v)}^-(t,k)\Big)$$ and we define
$$
 R_{osc,N}^{\varepsilon}=\sum_{k=0}^3
A_{k,N}^{\varepsilon}, $$ $$
R_{osc}^{\varepsilon,N}=R_{osc}^{\varepsilon}-R_{osc,N}^{\varepsilon},\quad
{\tilde R}_{osc,N}^{\varepsilon}=\sum_{k=0}^3 {\tilde
A}_{k,N}^{\varepsilon}. $$ Considering
${\varphi}_N^{\varepsilon}=W^{\varepsilon}+\varepsilon{\tilde
R}_{osc,N}^{\varepsilon}=({\varphi}^{\varepsilon}_{N,1},
{\varphi}^{\varepsilon}_{N,2})$, then $\varphi^{\varepsilon}_N$
satisfies the following equation $$
\partial_t\varphi^{\varepsilon}_N+Q^{\varepsilon}(\varphi^{\varepsilon}_N,
\varphi^{\varepsilon}_N+2\varepsilon\tilde{R}_{osc,N}^{\varepsilon}+2V)+
a_2^{\varepsilon}(D)\varphi^{\varepsilon}_N+{\tilde L}^{\varepsilon}
(\varphi^{\varepsilon}_N)=F^{\varepsilon} + R_{osc}^{\varepsilon,
N}+\varepsilon r_{osc, N}^{\varepsilon},
$$ where
$$\varepsilon r_{osc,
N}^{\varepsilon}=\varepsilon\Big(Q^{\varepsilon}
(\tilde{R}_{osc,N}^{\varepsilon},
\varepsilon\tilde{R}_{osc,N}^{\varepsilon}+2V)+a_2^{\varepsilon}(D)
\tilde{R}_{osc,N}^{\varepsilon}+{\tilde
L}^{\varepsilon}(\tilde{R}_{osc,N}^{\varepsilon})\Big)+\Big(R_{osc,N}^{\varepsilon}+\varepsilon\partial_t
\tilde{R}_{osc,N}^{\varepsilon}\Big).$$
Now, we have to studying the
low frequencies terms. This study is easy. In fact, we have the
following result.
\begin{lem}\label{l43}
A constant $C_N(T)$ exists, depending only $T$ and $N$ such that
$$ \|{\tilde R}^{\varepsilon}_{osc,N}\|_{{\mathcal C}^0_T(
H^{s-2})\cap L^2_T(H^{s-1})}\leq C_N(T), $$
 $$
\|r^{\varepsilon}_{osc,N}\|_{L^2_T(H^{s-2})}\leq C_N(T). $$
\end{lem}
\noindent{\bf Proof.} Let us recall that all the functions
considered here are truncated in low frequencies. Hence the result
is simply due to the fact that $v, b\in{\mathcal C}^0_T( H^s)$,
$\partial_t v,\partial_t b\in{\mathcal C}^0_T(H^{s-2})$ and the
following product law.\endproof
\begin{prop}\label{p235}
Let $s$ be an integer, $\sigma>\frac{3}{2}$. A constant $C$ exists
such that for all $f\in H^{\sigma}$ and $g\in H^{\sigma+1}$ with
$div(f)=div(g)=0$,
 we have $$ |{<
Q^{\varepsilon}(f,g),f>}_{H^{\sigma}}|\leq C\|f\|_{H^{\sigma}}^2
\|g\| _{H^{\sigma+1}} $$ and $$ |{<
Q^{\varepsilon}(f,f),f>}_{H^{\sigma}}|\leq C\|f\|_{H^{\sigma}}^3.
$$
\end{prop}
\noindent(The proof of the proposition used the fact ${\mathcal
L}(t)$ is a isometri in the Sobolev space and lemma \ref{l41}).
\begin{lem}\label{l44}
For any function $f\in {\mathcal C}^0_T( H^{s})$ with $s\in \R$, the
high frequency term $$ f^N={\mathcal F}^{-1}\Big({\bf
1}_{[N,+\infty[}{\mathcal F}(f)\Big) $$ goes to zero when $N$ goes
to infinity in ${\mathcal C}^0_T(H^{s}).$
\end{lem}
\noindent{\bf Proof.} Let us recall that $$
\|f^N(t)\|_{H^s}^2=\sum_{| k|\geq N}| k|^{2s}|{\mathcal
F}(f)(t,k)|^2. $$ So, the desired result is simply due to Dini's
theorem which implies that $\|f^N(t)\|_{H^s}^2$ goes to zero
uniformly in $t$.\endproof
\noindent This lemma implies in a
straightforward way the following result.
\begin{lem}\label{l45}
The high frequency term $R^{\varepsilon,N}_{osc}$ goes to zero in
${\mathcal C}^0_T( H^{s-2})\cap L^2_T(H^{s-1})$ when $N$ goes to
infinity, uniformly in $\varepsilon$; precisely
$\|R^{\varepsilon,N}_{osc}\|_{{\mathcal C}^0_T( H^{s-2})\cap
L^2_T(H^{s-1})}\leq \eta_{_N}$ with $\eta_{_N}\rightarrow 0$.
\end{lem}
\noindent Now, we can end the proof of the theorem. By the energy
estimate in $H^{s-2}({\mathbb{T}}^3)$ we obtain
\begin{eqnarray*}
\nonumber
\frac{1}{2}\frac{d}{dt}\|\varphi^{\varepsilon}_N\|_{H^{s-2}}^2+
\varepsilon\|\nabla\varphi^{\varepsilon}_{N,1}\|_{H^{s-2}}^2
+\sqrt{\varepsilon}\|\nabla\varphi^{\varepsilon}_{N,2}\|_{H^{s-2}}^2&\leq&
|{< F^{\varepsilon}+R_{osc}^{\varepsilon, N}+\varepsilon r_{osc,
N}^{\varepsilon},
\varphi^{\varepsilon}_N>}_{H^{s-2}}|\nonumber\\
&&+2|{< Q^{\varepsilon}(\varphi^{\varepsilon}_N,
\varepsilon\tilde{R}_{osc,N}^{\varepsilon}+V),
\varphi^{\varepsilon}_N>}_{H^{s-2}}|,\nonumber\\
&&+|{< Q^{\varepsilon}(\varphi^{\varepsilon}_N,
\varphi^{\varepsilon}_N),
\varphi^{\varepsilon}_N>}_{H^{s-2}}|\nonumber\\
\end{eqnarray*}
Which leads, using the product law (Proposition \ref{p235}), to
\begin{eqnarray}
\nonumber
\frac{1}{2}\frac{d}{dt}\|\varphi^{\varepsilon}_N\|_{H^{s-2}}^2\leq
C\Big[\|\varphi^{\varepsilon}_N\|_{H^{s-2}}^2(\|\varphi^{\varepsilon}_N\|_{H^{s-2}}+
\|\varepsilon\tilde{R}_{osc,N}^{\varepsilon}+V\|_{H^{s-1}}+1)\\
\nonumber +\|F^{\varepsilon} + R_{osc}^{\varepsilon, N}+\varepsilon
r_{osc, N}^{\varepsilon}\|_{H^{s-2}}^2\Big].
\end{eqnarray}
Integrating this inequality and using Lemma \ref{l43}, we obtain
\begin{eqnarray}
\nonumber &&\|\varphi^{\varepsilon}_N(t)\|_{H^{s-2}}^2
\leq\Big(\varepsilon C_N(T)^2+4\varepsilon^2\|\Delta
u\|_{L^2_{T}(H^{s-2})}^2+4\|R_{osc}^{\varepsilon,
N}\|_{L^2_{T}(H^{s-2})}^2\\ \nonumber &&+4\varepsilon^2
C_N(T)^2\Big)+C\int_0^t\|\varphi^{\varepsilon}_N(\tau)\|_{H^{s-2}}^2\Big(B(T)+\varepsilon
C_N(T)+\|\varphi^{\varepsilon}_N(\tau)\|_{H^{s-2}}\Big)d\tau,
\end{eqnarray}
where $B(T)=1+\|V\|_{L^{\infty}_{T}(H^{s-2})}.$\\ We set $$
T^*=\sup\{\quad 0\leq
t<T\quad/\quad\|\varphi^{\varepsilon}_N\|_{L_t^{\infty}(H^{s-2})}\leq
B(T)\quad\}. $$ Then, for all $0\leq t< T^*$, we can write,
\begin{eqnarray}
\nonumber \|\varphi^{\varepsilon}_N(t)\|_{H^{s-2}}^2&\leq
&\varepsilon C_N(T)^2+4\varepsilon^2\|u\|_{L^2_{T_0}(H^{s})}^2
+4\eta_{_N}^2+4\varepsilon^2 C_N(T)^2\\ \nonumber
&&\quad+c(B(T)+\varepsilon
C_N(T))\int_0^t\|\varphi^{\varepsilon}_N(\tau)\|_{H^{s-2}}^2d\tau.
\end{eqnarray}
A classical Gronwall estimate gives
\begin{eqnarray}
\nonumber
\|\varphi^{\varepsilon}_N(t)\|_{H^{s-2}}^2\leq\Big(\varepsilon
C_N(T)^2+4\varepsilon^2\|u\|_{L^2_{T}(H^{s})}^2
+4\eta_{_N}^2+4\varepsilon^2 C_N(T)^2\Big)\\ \nonumber \exp\Big(c T
B(T)+cT\varepsilon C_N(T)\Big).
\end{eqnarray}
For $N$ large and $\varepsilon$ small enough we obtain, thanks to
Lemma \ref{l45} $$
\|\varphi^{\varepsilon}_N(t)\|_{H^{s-2}}\leq\frac{B(T)}{2}, $$
 which easily implies that $T^*=T$. \\ Using Lemma \ref{l43}  and letting
 $\varepsilon\rightarrow 0$, $N\rightarrow+\infty$, we obtain
 $$
 W^{\varepsilon}\rightarrow 0\quad in\quad{\mathcal
 C}^0_T(H^{s-2}).
 $$
An interpolation argument concludes the proof of Theorem
\ref{t44}.\endproof
\section{The case $\Omega=\mathbb{R}^3$}
\noindent This section is devoted to study dispersion phenomena in
the ${\bf (MHD^\varepsilon)}$ system in the case of the space
$\R^3$. Let us introduce the ``linearized" equation in
$u^{\varepsilon}$ of the first equation of ${\bf
(MHD^\varepsilon)}$. $$ \left\{
\begin{array}{cllll}
&\partial_tu^{\varepsilon}+\displaystyle\frac{1}{\varepsilon}Lu^{\varepsilon}=-\nabla
p\quad in\quad \R_t\times \R_x^3\\
&\mbox{div}~u^{\varepsilon}=0\\ &u^{\varepsilon}(0)=u_0.
\end{array}
\right.$$ In Fourier variables $\xi\in \R^3$, we obtain $$
\partial_t{\mathcal F}(u^{\varepsilon})+
\frac{\xi_3}{\varepsilon|\xi|^2}\xi\times{\mathcal
F}(u^{\varepsilon})=0.$$ Hence, we are led to study the following
family of operators $$ {\mathcal
G}^{\varepsilon}:f\longmapsto\int_{\R_\xi^3}{\mathcal
F}(f)(\xi)\exp\Big({\mp it\frac{\xi_3}{\varepsilon|\xi|}+i
x.\xi}\Big) d\xi\quad\quad\quad\quad\quad\quad$$ $$
\quad\quad\quad\quad\quad\quad=\int_{\R_y^3\times
\R_\xi^3}f(y)\exp\Big({\mp
it\frac{\xi_3}{\varepsilon|\xi|}+i(x-y).\xi}\Big) d\xi dy. $$ We
notice that the phase function $\frac{\xi_3}{|\xi|}$ is almost
stationary when $\xi_3$ is almost equal to $0$ as well as when
$|\xi_3|$ is much larger then $|\xi_h|$. So, for some $0<~r~<~R$,
let us define the domain ${\mathcal C}_{r,R}$ by
$$ {\mathcal C}_{r,R}=\{\xi\in \R^3;\quad|\xi_3|>
r\quad\mbox{and}\quad|\xi|\leq R\} .$$ We consider $\psi$ a cut-off
function, which is radial with respect to horizontal variable
$\xi_h=(\xi_1,\xi_2)$ and whose value is 1 near
${\mathcal C}_{r,R}$.\\\\
First, we study the case when $\mathcal{F}(f)$ is supported in
${\mathcal C}_{r,R}$. We can write $$ {\mathcal
G}^{\varepsilon}f(t,x)=\Big(K(\frac{t}{\varepsilon},.)\ast
f\Big)(x),$$ where the kernel $K$ is defined by $$
K(t,z)=\int_{\R^3}\psi(\xi)e^{it\frac{\xi_3}{|\xi|}+iz.\xi}d\xi.
$$ As in \cite{CDGG}, we recall the following property of $K$.
\begin{lem}\label{l46}
For all $r$, $R$ such that $0< r< R$, there exists a constant
$C_{r,R}$ such that $$ \|K(t,.)\|_{L^{\infty}(\R^3)}\leq
C_{r,R}\mbox{min}\{1,t^{-\frac{1}{2}}\}. $$
\end{lem}
Let us denote by $w^{\varepsilon}$ the solution of $$ \left\{
\begin{array}{cllll}
&\partial_tw^{\varepsilon}+\displaystyle\frac{1}{\varepsilon}Lw^{\varepsilon}=f\quad\mbox{in}\quad
\R_t\times \R_x^3\\ &w^{\varepsilon}(0)=w_0.
\end{array}
\right. \leqno(PLF_{\varepsilon}) $$ Lemma \ref{l46} yields, in a
standard way, the following Strichartz-estimate (see \cite{CDGG}).
\begin{cor}\label{c242}
For all constants $r$ and $R$ such that $0< r< R$, let ${\mathcal
C}_{r,R}$ be the domain defined above. Then a constant $C_{r,R}$
exists such that if $$ \mbox{supp}{\mathcal
F}(w_0)\cup\mbox{supp}{\mathcal F}(f)\subset{\mathcal C}_{r,R} $$
then the solution $w^{\varepsilon}$ of $(PLF_{\varepsilon})$ with
the forcing term $f$ and initial data $w_0$ satisfies $$
\|w^{\varepsilon}\|_{L^4(\R_+,L^{\infty})}\leq
C_{r,R}\varepsilon^{\frac{1}{4}}\Big(\|w_0\|_{L^2}+\|f\|_{L^1(\R_+,L^2)}\Big).
$$
\end{cor}
\noindent We notice that the constant $C_{r,R}$ does not depend on
$\varepsilon$.\\\\
Using the above estimate, we are able to prove the convergence
result.
\subsection{Proof of Theorem \ref{t45}} We start by proving a convergence result for the
Leray's solutions. We shall apply the method used in \cite{CDGG} for
the first equation of our system.\\ Let $R> 0$ and $\chi$ be a
cut-off function in ${\mathcal D}(\R)$ taking the value $1$ near the
origin. We set
$u^{\varepsilon}_R:=\chi(\frac{|D|}{R})u^{\varepsilon}={\mathcal
F}^{-1}\Big(\chi(\frac{|\xi|}{R}){\mathcal
F}(u^{\varepsilon})(\xi)\Big)$. Then $u^{\varepsilon}_R$ satisfies
$$
\partial_t u^{\varepsilon}_R+
\frac{1}{\varepsilon} {\mathbb{P}}(u^{\varepsilon}_R\times
e_3)=f^{\varepsilon}_R, $$ with $$ f^{\varepsilon}_R=-\chi(\frac{|
D|}{R})\Big(\varepsilon\Delta
u^{\varepsilon}+{\mathbb{P}}(u^{\varepsilon}\nabla
u^{\varepsilon})-{\mathbb{P}}(b^{\varepsilon}\nabla
b^{\varepsilon})-\sqrt{\varepsilon}\partial_3 b^{\varepsilon}\Big).
$$ By Duhamel's formula, we can write $$
u^{\varepsilon}_R=e^{t\sigma^{\varepsilon}(D)}u_{0,R}+\int_0^t
e^{(t-t')\sigma^{\varepsilon}(D)}f^{\varepsilon}_R(t')d t',$$ where
$\sigma^{\varepsilon}(\xi)u=\frac{\xi_3}{\varepsilon|\xi|^2}u\times\xi$.\\
So we have
\begin{eqnarray*}
\|u^{\varepsilon}_R\|_{L^4_T(L^{\infty}(R^3))}\leq
\|e^{t\sigma^{\varepsilon}(D)}u_{0,R}\|_{L^4_T(L^{\infty}(R^3))}+
\|\int_0^te^{(t-t')\sigma^{\varepsilon}(D)}f^{\varepsilon}_R(t')dt'\|_{L^4_T(L^{\infty})}.
\end{eqnarray*} For the first term of the right side of the above inequality,
we localize on ${\mathcal C}_{r,R}$. Hence, for all $0< r< R$, we
write
\begin{eqnarray}
\|e^{t\sigma^{\varepsilon}(D)}u_{0,R}\|_{L^4_T(L^{\infty}(R^3))}\leq
\|e^{t\sigma^{\varepsilon}(D)}\chi(\frac{D_3}{r})u_{0,R}\|_{L^4_T(L^{\infty}(R^3))}
+\|e^{t\sigma^{\varepsilon}(D)}(Id-\chi(\frac{D_3}{r})
)u_{0,R}\|_{L^4_T(L^{\infty}(R^3))}. \nonumber\end{eqnarray} Using
Corollary \ref{c242} we obtain

\begin{equation}\label{S1}
\|e^{t\sigma^{\varepsilon}(D)}u_{0,R}\|_{L^4_T(L^{\infty}(R^3))}\leq
CT^{\frac{1}{4}}(R^2r)^{\frac{1}{2}}\|u_0\|_{L^2}+
{\varepsilon}^{\frac{1}{4}}C_{r,R}\|u_0\|_{L^2}.
\end{equation}
For the other term we have
\begin{eqnarray}
\|\int_0^te^{(t-t')\sigma^{\varepsilon}(D)}f^{\varepsilon}_R
(t')dt'\|_{L^4_T(L^{\infty})}&\leq  &
\|\int_0^te^{(t-t')\sigma^{\varepsilon}(D)}\chi(\frac{D_3}{r})
f^{\varepsilon}_R(t')dt'\|_{L^4_T(L^{\infty})}
\nonumber\\&&+\|\int_0^te^{(t-t')
\sigma^{\varepsilon}(D)}(Id-\chi(\frac{D_3}{r}))f^{\varepsilon}_R(t')dt'\|_{L^4_T(L^{\infty})}.
\nonumber\end{eqnarray} Corollary \ref{c242} and Bernstein's lemma
imply the following
\begin{eqnarray}
\|\int_0^te^{(t-t')\sigma^{\varepsilon}(D)}f^{\varepsilon}_R(t')dt'\|_{L^4_T(L^{\infty})}&\leq
& CT^{\frac{5}{4}}R^{\frac{3}{2}}\|\chi(\frac{D_3}{r})\chi(\frac{|
D|}{R})f^{\varepsilon} \|_{L^{\infty}_T(L^2)}\nonumber\\
&&+C_{r,R}\varepsilon^{\frac{1}{4}}\|(Id-\chi(\frac{D_3}{r}))\chi(\frac{|
D|}{R})f^{\varepsilon} \|_{L^1_T(L^2)}\nonumber\\ &\leq  &\!
CT^{\frac{5}{4}}R^{\frac{3}{2}}\|\chi(\frac{D_3}{r})\chi(\frac{|
D|}{R})f^{\varepsilon}
\|_{L^{\infty}_T(L^2)}\nonumber\\&&+C_{r,R}\varepsilon^{\frac{1}{4}}\|\chi(\frac{|
D|}{R})f^{\varepsilon} \|_{L^{\infty}_T(L^2)}. \nonumber
\end{eqnarray}
Now, we focus our attention to the first right term of the last
inequality. \begin{eqnarray*} \chi(\frac{D_3}{r})\chi(\frac{|
D|}{R})f^{\varepsilon}=&&\varepsilon\chi(\frac{D_3}{r})\chi(\frac{|
D|}{R})\Delta u^\varepsilon-\chi(\frac{D_3}{r})\chi(\frac{|
D|}{R}){\mathbb{P}}(u^{\varepsilon}.\nabla
u^{\varepsilon})\\&+&\chi(\frac{D_3}{r})\chi(\frac{|
D|}{R}){\mathbb{P}}(b^{\varepsilon}.\nabla
b^{\varepsilon})+\sqrt{\varepsilon}\chi(\frac{D_3}{r})\chi(\frac{|
D|}{R})\partial_3 b^{\varepsilon}. \end{eqnarray*} First, notice
that we have
\begin{eqnarray} \|\chi(\frac{D_3}{r})\chi(\frac{|
D|}{R}){\mathbb{P}}(u^{\varepsilon}.\nabla
u^{\varepsilon})\|_{L^{\infty}_T(L^2)}&\leq  &
R\|\chi(\frac{D_3}{r})\chi(\frac{| D|}{R})(u^{\varepsilon}\otimes
u^{\varepsilon})\|_{L^{\infty}_T(L^2)}\nonumber\\&\leq  &
 R\|{\mathcal
F}^{-1}\chi(\frac{\xi_3}{r})\chi(\frac{|
\xi|}{R})*(u^{\varepsilon}\otimes
u^{\varepsilon})\|_{L^{\infty}_T(L^2)}\nonumber\\&\leq  &
R\|{\mathcal F}^{-1}\chi(\frac{\xi_3}{r})\chi(\frac{|
\xi|}{R})\|_{L^{\infty}_T(L^2)}\|u^{\varepsilon}\otimes
u^{\varepsilon}\|_{L^{\infty}_T(L^1)}\nonumber\\&\leq  &
R(R^2r)^{\frac{1}{2}}\|U_0\|_{L^2}^2. \nonumber\end{eqnarray} Then,
in a same manner as above, we obtain $$
\|\chi(\frac{D_3}{r})\chi(\frac{|
D|}{R}){\mathbb{P}}(b^{\varepsilon}.\nabla
b^{\varepsilon})\|_{L^{\infty}_T(L^2)}\leq
R(R^2r)^{\frac{1}{2}}\|U_0\|_{L^2}^2 $$ and finally $$
\|\chi(\frac{D_3}{r})\chi(\frac{| D|}{R})\partial_3
b^{\varepsilon}\|_{L^{\infty}_T(L^2)}\leq
r\|\chi(\frac{D_3}{r})\chi(\frac{|
D|}{R})b^{\varepsilon}\|_{L^{\infty}_T(L^2)}\leq
r\|b^{\varepsilon}\|_{L^{\infty}_T(L^2)}\leq r\|U_0\|_{L^2}^2.
$$ For the second right term we use
\begin{eqnarray}
\|\chi(\frac{| D|}{R})f^{\varepsilon}\|_{L^{\infty}_T(L^2)}&\leq
&\|\chi(\frac{| D|}{R}){\mathbb{P}}(u^{\varepsilon}.\nabla
u^{\varepsilon})\|_{L^{\infty}_T(L^2)}+\|\chi(\frac{|
D|}{R}){\mathbb{P}}(b^{\varepsilon}.\nabla
b^{\varepsilon})\|_{L^{\infty}_T(L^2)}\nonumber\\&&+\varepsilon\|\chi(\frac{|
D|}{R})\Delta u
^{\varepsilon}\|_{L^{\infty}_T(L^2)}+\sqrt{\varepsilon}\|\chi(\frac{|
D|}{R})\partial_3
b^{\varepsilon}\|_{L^{\infty}_T(L^2)}\nonumber\\&\leq  & C
R^{\frac{5}{2}}\|U_0\|_{L^2}^2 +C
R^{\frac{5}{2}}\|U_0\|_{L^2}^2+R\|b^{\varepsilon}\|_{L^{\infty}_T(L^2)}\nonumber\\&\leq
& C R^{\frac{5}{2}}\|U_0\|_{L^2}^2+R\|U_0\|_{L^2}.
\nonumber\end{eqnarray} Finally, we have
\begin{eqnarray}
\|\int_0^t
e^{{(t-t')}(\sigma^{\varepsilon}(D)}f^{\varepsilon}_R(t')dt'\|_{L^4_T(L^{\infty})}&\leq
& \!CT^{\frac{5}{4}}R^{\frac{3}{2}}(\|U_0\|_{L^2}^2 R^2
r^{\frac{1}{2}}+r\|U_0\|_{L^2})\nonumber\\
&&+C_{r,R}\varepsilon^{\frac{1}{4}}T
(R^{\frac{5}{2}}\|U_0\|_{L^2}^2+R\|U_0\|_{L^2}).
\label{S2}\end{eqnarray} Inequalities (\ref{S1}) and (\ref{S2})
imply in a simple manner that
\begin{eqnarray}
u_R^{\varepsilon}\longrightarrow0\quad \mbox{as}\quad
\varepsilon\rightarrow0\; ;\; in \quad L^4_T(L^{\infty}). \label{S3}
\end{eqnarray}Using the embedding $
L^{\infty}\hookrightarrow{\mathcal C}^{-\frac{1}{4}}$, we deduce
that $$u_R^{\varepsilon}\longrightarrow0\quad as\quad
\varepsilon\rightarrow0\; ;\;  in \quad L^4_T({\mathcal
C}^{-\frac{1}{4}}).$$ On the other hand, we have
\begin{eqnarray}
\|u^{\varepsilon}-u_R^{\varepsilon}\|_{L^2_T({\mathcal
C}^{-\frac{1}{4}})}&\leq  &
C\|u^{\varepsilon}-u_R^{\varepsilon}\|_{L^2_T(H^{-\frac{7}{4}})}\nonumber\\
&\leq &CR^{-\frac{7}{4}}\|u^{\varepsilon}\|_{L^2_T(L^2(\R^3))}
\nonumber
\end{eqnarray}
and then by energy estimate \ref{eef}, we conclude that
\begin{eqnarray}
\|u^{\varepsilon}-u_R^{\varepsilon}\|_{L^2_T({\mathcal
C}^{-\frac{1}{4}})}&\leq  &
CR^{-\frac{7}{4}}T^{\frac{1}{2}}\|U_0\|_{L^2(\R^3)}. \label{S4}
\end{eqnarray}
Hence, inequalities (\ref{S3}) and (\ref{S4}) and an interpolation
argument, we can deduce the first part of Theorem \ref{t45}. \\
Now, we come to the proof of the second part of Theorem \ref{t45}.
We recall that we have $$
\partial_t b^{\varepsilon}=\varepsilon\Delta b^{\varepsilon}+b^{\varepsilon}.\nabla
u^{\varepsilon}-u^{\varepsilon}.\nabla
b^{\varepsilon}-\sqrt{\varepsilon}\partial_3 u^{\varepsilon}. $$ We
take the scalar produce in $H^{s-3}$ with $b^\varepsilon-b_0$, we
obtain
\begin{eqnarray*}
\|b^{\varepsilon}(t)-b_0\|_{H^{s-3}}^2&\leq  &\int_0^t|{<
b^{\varepsilon}(\tau).\nabla
u^{\varepsilon}(\tau),b^{\varepsilon}(\tau)-b_0>}_{H^{s-3}}|d\tau\nonumber
\\
&&+\int_0^t|{< u^{\varepsilon}(\tau).\nabla
b^{\varepsilon}(\tau),b^{\varepsilon}(\tau)-b_0>}_{H^{s-3}}|d\tau
\nonumber
\\
&&+\varepsilon\int_0^t\|\Delta
b^{\varepsilon}(\tau)\|_{H^{s-3}}\|b^{\varepsilon}(\tau)-b_0\|_{H^{s-3}}d\tau\nonumber
\\
&&+\sqrt{\varepsilon}\int_0^t\|\partial_3b^{\varepsilon}(\tau)\|_{H^{s-3}}
\|b^{\varepsilon}(\tau)-b_0\|_{H^{s-3}}d\tau,\nonumber
\\
\end{eqnarray*}
using the energy estimate (\ref{eef}) and the following product law
\begin{lem}\label{l47}
Let $\sigma>4$ an integer and $a$, $b$ and $c$ three vectors field
such that $a\in H^{\sigma}(\mathbb {R}^3)\cap {\mathcal
C}^{\sigma+\frac{1}{4}}(\mathbb {R}^3)$, $b\in H^{\sigma+1}(\mathbb
{R}^3)\cap {\mathcal C}^{\sigma+\frac{5}{4}}(\mathbb {R}^3)$, $c\in
H^{\sigma}(\mathbb {R}^3)$ and $\mbox{div}~a=0$. Then, a constant
$C(\sigma)$ exists such that $$ |{< a.\nabla b,c>}_{H^\sigma}|\leq
C(\sigma)\|c\|_{H^\sigma}\min\Big(\|a\|_{{\mathcal
C}^{\sigma+\frac{1}{4}}}\|\nabla b \|_{H^\sigma};\|a
\|_{H^\sigma}\|b\|_{{\mathcal C}^{\sigma+\frac{5}{4}}}\Big).$$
\end{lem}
 We obtain
\begin{eqnarray*}
\|b^{\varepsilon}(t)-b_0\|_{H^{s-3}}^2&\leq
&C\sqrt{T}\|(w_0,b_0)\|_{H^{s}}^2\|u^\varepsilon\|_{L_T^2({\mathcal
C}^{s-{\frac{7}{4}}})} \nonumber
\\
&&+C\sqrt{\varepsilon}\|(w_0,b_0)\|_{H^{s}}^2,\nonumber
\end{eqnarray*}
the fact $$u^\varepsilon\rightarrow 0\quad in\quad L^4_T({\mathcal
C}^{s-\frac{7}{4}})$$ and an interpolation argument, we can deduce
the desired result.\endproof
\subsection{Proof of Theorem \ref{t47}} We used the
following lemma
\begin{lem}\label{l48}For$\quad\bar{a}\in \big(H^{\sigma}({\mathbb R}^2)\big)^3$ and
$b,c\in\big(H^\sigma({\mathbb R}^3)\big)^3$ ($\sigma>eq 4$ an
integer) such that $\mbox{div}_h~\bar{a}=0$, then
$$
\begin{array}{cllll}
|{<{\bar a}.\nabla b ,c>}_{H^\sigma}+{<{\bar a}.\nabla c
,b>}_{H^\sigma}|&\leq  & C\|{\bar a}\|_{H^\sigma}\|b
\|_{H^{\sigma}}\|c \|_{H^{\sigma}}\\
|{<{\bar{a}}.\nabla b,b>}_{H^\sigma}|&\leq  &
C\|{\bar{a}}\|_{H^{\sigma}}\|b \|_{H^{\sigma}}^2\\
\end{array}
$$
Moreover, if $\bar{a}\in \big(H^{\sigma+1}({\mathbb R}^2)\big)^3$,
then$$|{< b.\nabla {\bar{a}}, c>}_{H^\sigma}|\leq C\|{\bar
a}\|_{H^{\sigma+1}}\|b \|_{H^{\sigma}}\|c \|_{H^{\sigma}}$$
\end{lem}
and the approximate system
$$ \left\{
\begin{array}{cllll}
\partial _tw -\varepsilon\Delta J_nw
+\displaystyle\frac{1}{\varepsilon}{J_nw\times e_3
}+\sqrt{\varepsilon}\partial_3J_nB-J_n(J_nB.\nabla J_nB)&&\\
\quad\quad\quad\quad\quad\quad=\nabla
\Delta^{-1}\mbox{div}~\Big(J_n(J_nB.\nabla
J_nB)-\displaystyle\frac{1}{\varepsilon}{J_nw\times e_3 }\Big)&&\\
\partial_tB^{\varepsilon}-\displaystyle{\sqrt{\varepsilon}}\Delta J_nB
+\sqrt{\varepsilon}\partial_3J_nw+J_n((J_n{\bar
u}^\varepsilon+J_nw).\nabla J_nB)&in&
[0,T_0]\times\mathbb{R}^3\\
\quad\quad\quad\quad\quad\quad\quad\quad\quad\quad\quad
-J_n(J_nB.\nabla (J_n{\bar u }^\varepsilon+J_nw))=0&&\\
(w,B)(0)=(J_nw_0,J_nb_0)&&\\
(\mbox{div}~w,\mbox{div}~B)=(0,0).&&\\
\end{array}
\right. \leqno({\bf MHD3D_n})$$ The Lemmas \ref{l41} and \ref{l42}
implies that $({\bf MHD3D}_n)$ have a ordinary differential equation
form in $H^s({\mathbb R}^3)$, then there exists $W_n:=(w_n,B_n)\in
{\mathcal C}^1([0,T^*_{n,\varepsilon}),H^s({\mathbb R}^3))$ a
maximal solution of $({\bf MHD3D}_n)$. The uniqueness and the fact
$J_n^2=J_n$ implies that
$J_nW_n=W_n$.\\
 Moreover
$T^*_{n,\varepsilon}> T_1:=\displaystyle\frac{1}{C(s)(\|{\bar
u}_0\|_{H^{s+1}}+\|(w_0,b_0)\|_{H^s})}$ and for all $t\in [0,T_1]$

\begin{equation}
\label{e42} \|W_n(t)\|^2_{H^s}+2\varepsilon\int_0^t\|\nabla
w_n\|^2_{H^s} +2\sqrt\varepsilon\int_0^t\|\nabla B_n\|^2_{H^s}\leq
2\|(w_0,b_0)\|^2_{H^s}.
\end{equation}
In the end we used the proof of theorem \ref{t45}, we obtain
existence of a $W^\varepsilon:=(w^\varepsilon,B^\varepsilon)\in
L^\infty_{T_1}(H^s)\cap L^2_{T_1}(H^{s+1})$ solution of $({\bf
{MHD3D}^\varepsilon})$, satisfy for all $t\in [0,T_1]$
\begin{equation} \label{e43}
\|W^\varepsilon(t)\|^2_{H^s}+2\varepsilon\int_0^t\|\nabla
w^\varepsilon\|^2_{H^s}+2\sqrt\varepsilon\int_0^t\|\nabla
B^\varepsilon\|^2_{H^s}\leq 2\|(w_0,b_0)\|^2_{H^s}.
\end{equation}
The proof of the uniqueness is easy.\endproof
\subsection{Proof of Theorem \ref{t48}}
We posed $v^\varepsilon:=u^\varepsilon-{\bar u}^\varepsilon$, then
$(v^\varepsilon,b^\varepsilon)$ is solution of the following system
$$\left\{
\begin{array}{cllll}
\partial_tv^{\varepsilon}-\varepsilon\Delta v^{\varepsilon}+
v^{\varepsilon}.\nabla v^{\varepsilon}+ v^{\varepsilon}.\nabla {\bar
u}^{\varepsilon}+ {\bar u}^{\varepsilon}.\nabla
v^{\varepsilon}\quad\quad\quad&&\\ \quad\quad\quad -\mbox{curl }
b^{\varepsilon}\times b^{\varepsilon}+\sqrt{\varepsilon}\mbox{curl
}b^{\varepsilon}\times e_3+\displaystyle\frac{v^{\varepsilon}\times
e_3}{\varepsilon}&=&-\nabla p^{\varepsilon}\\
\partial _t b^{\varepsilon} - \sqrt{\varepsilon}\Delta b^{\varepsilon}+
{\bar u}^{\varepsilon}. \nabla b^{\varepsilon}+v^{\varepsilon}.
\nabla b^{\varepsilon}-b^{\varepsilon}. \nabla
v^{\varepsilon}\quad\quad\quad &&\\
\quad\quad\quad\quad\quad\quad-b^{\varepsilon}.\nabla{\bar
u}^{\varepsilon} +\sqrt{\varepsilon}\mbox{curl
}(v^{\varepsilon}\times e_3)&=&0\\
\mbox{div}~v^{\varepsilon}&=&0\\
\mbox{div}~b^{\varepsilon}&=&0\\
(v^\varepsilon,b^\varepsilon)_{/{t=0}}&=&(w_0,b_0)\\
\end{array}\right.\leqno{({\bf {MHDF}^\varepsilon})}$$
Now we consider the approximate system $$\left\{
\begin{array}{cllll}
\partial_tv&-&\varepsilon\Delta J_n v+
J_n(J_nv.\nabla J_nv)+ J_n(J_nv.\nabla J_n{\bar u}^{\varepsilon}) +
J_n(J_n{\bar u}^{\varepsilon}.\nabla J_nv)\\&+&J_n(J_n{\bar
u}^{\varepsilon}.\nabla J_nv)-J_n(\mbox{curl} J_nb\times J_nb)
+\sqrt{\varepsilon}\mbox{curl }J_nb\times e_3+\\
\displaystyle\frac{J_nv\times e_3}{\varepsilon}&=&-\nabla
\Delta^{-1} \mbox{div}~\Big(J_n(J_nv.\nabla J_nv)+ J_n(J_nv.\nabla
J_n{\bar u}^{\varepsilon})+ J_n(J_n{\bar u}^{\varepsilon}.\nabla
J_nv)\\+\displaystyle\frac{J_nv\times
e_3}{\varepsilon}&+&J_n(J_n{\bar u}^{\varepsilon}.\nabla
J_nv)-J_n(\mbox{curl} J_nb\times J_nb)+\sqrt{\varepsilon}\mbox{curl
}J_nb\times e_3\Big),
\\
\partial _t b &-& \sqrt{\varepsilon}\Delta J_nb+
J_n(J_n{\bar u}^{\varepsilon}. \nabla J_nb)+J_n(J_nv. \nabla
J_nb)\quad in\quad [0,T_0]\times\mathbb{R}^3\\ &-&J_n(J_nb. \nabla
J_nv)-J_n(J_nb. \nabla J_nv)-J_n(J_nb.\nabla J_n{\bar
u}^{\varepsilon})\\&+&\sqrt{\varepsilon}\mbox{curl }(J_nv\times
e_3)=0\\ \mbox{div}~v&=&\mbox{div}~b=0\\
(v,b)(0)&=&(J_nw_0,J_nb_0)\\
\end{array}\right.\leqno{({\bf MHDF}_n^\varepsilon)}$$
We use Lemmas \ref{l41}, \ref{l43} we proof that $({\bf
MHDF}_n^\varepsilon)$ have a ordinary differential equation form in
$H^s({\mathbb R}^3)$, then there exists $V_n:=(v_n,b_n)\in {\mathcal
C}^1([0,T^*_{n,\varepsilon}),H^s({\mathbb R}^3))$ a maximal solution
of $({\bf MHDF}_n^\varepsilon)$. The uniqueness and the fact
$J_n^2=J_n$ implies that $J_nV_n=V_n$. The Lemmas \ref{l41} and
\ref{l42} implies $T^*_{n,\varepsilon}> T_1$ and for all $t\in
[0,T_1]$
\begin{equation}\label{e42}
\|V_n(t)\|^2_{H^s}+2\varepsilon\int_0^t\|\nabla v_n\|^2_{H^s}
+2\sqrt\varepsilon\int_0^t\|\nabla b_n\|^2_{H^s}\leq
2\|(w_0,b_0)\|^2_{H^s}.
\end{equation}
The end is similar of the proof for Theorem \ref{t45}, then we
obtain existence of a
$V^\varepsilon:=(v^\varepsilon,b^\varepsilon)\in
L^\infty_{T_1}(H^s)\cap L^2_{T_1}(H^{s+1})$ solution of $({\bf
{MHDF}^\varepsilon})$, satisfy for all $t\in [0,T_1]$
\begin{equation}\label{e43}
\|V^\varepsilon(t)\|^2_{H^s}+2\varepsilon\int_0^t\|\nabla
v^\varepsilon\|^2_{H^s} +2\sqrt\varepsilon\int_0^t\|\nabla
b^\varepsilon\|^2_{H^s}\leq 2\|(w_0,b_0)\|^2_{H^s}.
\end{equation}
The proof of the uniqueness is easy.\endproof
\subsection{Proof of Theorem \ref{t49}}
We posed $\psi^\varepsilon:=(u^\varepsilon-{\bar
u}^\varepsilon-w^\varepsilon,b^\varepsilon-B^\varepsilon)
=(\psi^\varepsilon_1,\psi^\varepsilon_2)$, $\psi^\varepsilon$
satisfy the followings equations \begin{eqnarray*}\partial_t
\psi_1^\varepsilon-\varepsilon\Delta
\psi_1^\varepsilon+\sqrt\varepsilon\partial_3\psi_2^\varepsilon&=&f^\varepsilon_1\\
\partial_t
\psi_2^\varepsilon-\sqrt\varepsilon\Delta
\psi_2^\varepsilon+\sqrt\varepsilon\partial_3\psi_1^\varepsilon&=&f^\varepsilon_2,\\
\end{eqnarray*}
where
$$\quad\quad f_1^\varepsilon:=-\nabla
p^\varepsilon-\psi_1^\varepsilon.\nabla
\psi_1^\varepsilon-\psi_1^\varepsilon.\nabla {\bar
u}^\varepsilon-{\bar u}^\varepsilon.\nabla
\psi_1^\varepsilon+\psi_2^\varepsilon.\nabla
\psi_2^\varepsilon+\psi_2^\varepsilon.\nabla
B^\varepsilon+B^\varepsilon.\nabla
\psi_1^\varepsilon-w^\varepsilon.\nabla {\bar u}^\varepsilon-{\bar
u}^\varepsilon.\nabla w^\varepsilon
$$
$$f_2^\varepsilon:=-\psi_1^\varepsilon.\nabla
\psi_2^\varepsilon-w^\varepsilon.\nabla \psi_2^\varepsilon-{\bar
u}^\varepsilon.\nabla \psi_2^\varepsilon-\psi_1^\varepsilon.\nabla
B^\varepsilon+B^\varepsilon.\nabla
\psi_1^\varepsilon+\psi_2^\varepsilon.\nabla {\bar
u}^\varepsilon+\psi_2^\varepsilon.\nabla
\psi_1^\varepsilon+\psi_2^\varepsilon.\nabla w^\varepsilon.
$$
We take the scalar produce in $H^{s-3}({\mathbb R}^3)$ we obtain
\begin{equation}\label{ef}\frac{1}{2}\frac{d}{dt}\|\psi^\varepsilon\|^2_{H^{s-3}}
+\varepsilon\|\nabla\psi_1^\varepsilon\|^2_{H^{s-3}}
+\sqrt\varepsilon\|\nabla\psi_2^\varepsilon\|^2_{H^{s-3}}\leq
\sum_{k=1}^{14}I_k,\end{equation} where \begin{eqnarray*}
&&I_1:=|{<\psi_1^\varepsilon.\nabla
\psi_1^\varepsilon,\psi_1^\varepsilon>}_{H^{s-3}}|,\quad
I_2:=|{<{\bar u}^\varepsilon.\nabla
\psi_1^\varepsilon,\psi_1^\varepsilon>}_{H^{s-3}}|\\&&I_3:=|{<\psi_1^\varepsilon.\nabla
{\bar u }^\varepsilon,\psi_1^\varepsilon>}_{H^{s-3}}|,\quad
I_4:=|{<\psi_2^\varepsilon.\nabla
B^\varepsilon,\psi_1^\varepsilon>}_{H^{s-3}}|\\&&I_5:=|{<\psi_2^\varepsilon.\nabla
\psi_2^\varepsilon,\psi_1^\varepsilon>}_{H^{s-3}}+{<\psi_2^\varepsilon.\nabla
\psi_1^\varepsilon,\psi_2^\varepsilon>}_{H^{s-3}}|\\&&I_6:=|{<
B^\varepsilon.\nabla
\psi_2^\varepsilon,\psi_1^\varepsilon>}_{H^{s-3}}+{<
B^\varepsilon.\nabla
\psi_1^\varepsilon,\psi_2^\varepsilon>}_{H^{s-3}}|\\
&&I_7:=|{<\psi_1^\varepsilon.\nabla
\psi_2^\varepsilon,\psi_2^\varepsilon>}_{H^{s-3}}|,\quad
I_8:=|{<\psi_1^\varepsilon.\nabla
B^\varepsilon,\psi_2^\varepsilon>}_{H^{s-3}}|\\ && I_9:=|{<{\bar u
}^\varepsilon.\nabla
\psi_2^\varepsilon,\psi_2^\varepsilon>}_{H^{s-3}}|,\quad I_{10}:=|{<
w^\varepsilon.\nabla
\psi_2^\varepsilon,\psi_2^\varepsilon>}_{H^{s-3}}|\\
&&I_{11}:=|{<\psi_2^\varepsilon.\nabla {\bar u
}^\varepsilon,\psi_2^\varepsilon>}_{H^{s-3}}|,\quad
I_{12}:=|{<\psi_2^\varepsilon.\nabla
w^\varepsilon,\psi_1^\varepsilon>}_{H^{s-3}}|\\ && I_{13}:=|{<{\bar
u }^\varepsilon.\nabla
w^\varepsilon,\psi_1^\varepsilon>}_{H^{s-3}}|,\quad I_{14}:=|{<
w^\varepsilon.\nabla {\bar u
}^\varepsilon,\psi_1^\varepsilon>}_{H^{s-3}}|.\\
\end{eqnarray*}
Using Cauchy-Schwarz inequality and the fact $H^{s-3}$ is an algebra
we obtain
\begin{eqnarray}
I_4&\leq  &C\|B^\varepsilon\|_{H^{s-2}}
\|\psi_2^\varepsilon\|_{H^{s-3}}
\|\psi_1^\varepsilon\|_{H^{s-3}}\leq
C\|\psi^\varepsilon\|_{H^{s-3}}^2
\\ I_8&\leq  &C\|B^\varepsilon\|_{H^{s-2}}
\|\psi_1^\varepsilon\|_{H^{s-3}}
\|\psi_2^\varepsilon\|_{H^{s-3}}\leq
C\|\psi^\varepsilon\|_{H^{s-3}}^2\\
I_{12}&\leq  &C\|w^\varepsilon\|_{H^{s-2}}
\|\psi_1^\varepsilon\|_{H^{s-3}}
\|\psi_2^\varepsilon\|_{H^{s-3}}\leq
C\|\psi^\varepsilon\|_{H^{s-3}}^2.\\ \nonumber
\end{eqnarray}
Using lemma \ref{l48} we obtain
\begin{eqnarray}
I_1&\leq  &C\|\psi_1^\varepsilon\|_{H^{s-3}}^3\leq
C\|\psi^\varepsilon\|_{H^{s-3}}^2
\\
I_2&\leq  &C\|{\bar u}^\varepsilon\|_{H^{s-2}}
\|\psi_1^\varepsilon\|_{H^{s-3}}^2\leq
C\|\psi^\varepsilon\|_{H^{s-3}}^2
\\I_3&\leq  &C\|{\bar u}^\varepsilon\|_{H^{s-2}}
\|\psi_1^\varepsilon\|_{H^{s-3}}^2\leq
C\|\psi^\varepsilon\|_{H^{s-3}}^2
\\I_5&\leq  &C\|\psi_1^\varepsilon\|_{H^{s-2}}
\|\psi_2^\varepsilon\|_{H^{s-3}}^2\leq
C\|\psi^\varepsilon\|_{H^{s-3}}^2
\\I_6&\leq  &C\|B^\varepsilon\|_{H^{s-3}}
\|\psi_2^\varepsilon\|_{H^{s-3}}
\|\psi_1^\varepsilon\|_{H^{s-3}}\leq
C\|\psi^\varepsilon\|_{H^{s-3}}^2
\\I_7&\leq  &C\|\psi_1^\varepsilon\|_{H^{s-3}}
\|\psi_2^\varepsilon\|_{H^{s-3}}^2\leq
C\|\psi^\varepsilon\|_{H^{s-3}}^2
\\I_9&\leq  &C\|{\bar u}^\varepsilon\|_{H^{s-3}}
\|\psi_2^\varepsilon\|_{H^{s-3}}^2\leq
C\|\psi^\varepsilon\|_{H^{s-3}}^2
\\I_{10}&\leq  &C\|w^\varepsilon\|_{H^{s-2}}
\|\psi_2^\varepsilon\|_{H^{s-3}}^2\leq
C\|\psi^\varepsilon\|_{H^{s-3}}^2
\\I_{11}&\leq  &C\|{\bar u}^\varepsilon\|_{H^{s-2}}
\|\psi_2^\varepsilon\|_{H^{s-3}}^2\leq
C\|\psi^\varepsilon\|_{H^{s-3}}^2.
\\
\nonumber
\end{eqnarray}
The difficult is the terms $I_{13}$ and $I_{14}$, for exemple study
the term $I_{14}$.\\ Set $\chi\in{\mathcal D}({\mathbb R})$ valu $1$
near the origin. For $R> 0$ we pose
$$w^\varepsilon_R:={\mathcal
F}^{-1}(\chi(\frac{|\xi|}{R}){\mathcal
F}(w^\varepsilon)),\quad{\tilde
w}^\varepsilon_R:=w^\varepsilon-w^\varepsilon_R.$$
 We have
$$I_{14}\leq J+{\tilde J},$$ where $J:=|{< w^\varepsilon_R.\nabla {\bar u
}^\varepsilon;\psi_1^\varepsilon>}_{H^{s-3}}|$, ${\tilde
J}:=|{<{\tilde w}_R^\varepsilon.\nabla {\bar u
}^\varepsilon;\psi_1^\varepsilon>}_{H^{s-3}}|.$\\  $\bullet$ Study
of $J$
\begin{eqnarray*}
J&=&|\sum_{|\alpha|\leq{s-3}}\int
 \partial^\alpha( w_R^\varepsilon.\nabla {\bar u}^\varepsilon)\partial^\alpha
 \psi_1^\varepsilon|\\
 &=&|\sum_{|\alpha|\leq{s-3}}\sum_{\beta\leq\alpha}C^\beta_\alpha\int
 \partial^\beta w_R^\varepsilon.\nabla \partial
 ^{\alpha-\beta}{\bar u}^\varepsilon\partial^\alpha
 \psi_1^\varepsilon|\\
 &\leq  &\sum_{|\alpha|\leq{s-3}}\sum_{\beta\leq\alpha}C^\beta_\alpha|\int
 \partial^\beta w_R^\varepsilon.\nabla \partial
 ^{\alpha-\beta}{\bar u}^\varepsilon\partial^\alpha \psi_1^\varepsilon|.
\end{eqnarray*}
For any $\alpha, \beta$ we have
\begin{eqnarray*}
|\int
 \partial^\beta w_R^\varepsilon.\nabla \partial
 ^{\alpha-\beta}{\bar u}^\varepsilon\partial^\alpha
 \psi_1^\varepsilon|&=&
|\int
 {\mathcal F}(\partial^\beta w_R^\varepsilon)\chi(\frac{|\xi|}{2R}){\overline{
 {\mathcal F}(\nabla \partial
 ^{\alpha-\beta}{\bar u}^\varepsilon\partial^\alpha
 \psi_1^\varepsilon)}}|\\
 &=&|\int
 (\partial^\beta w_R^\varepsilon){\mathcal F}^{-1}\Big(\chi(\frac{|\xi|}{2R})
 {\overline{
 {\mathcal F}(\nabla \partial
 ^{\alpha-\beta}{\bar u}^\varepsilon\partial^\alpha
 \psi_1^\varepsilon)}}\Big)|\\&\leq  &\|\partial^\beta w_R^\varepsilon\|_{L^\infty}
 \|{\mathcal F}^{-1}(\chi(\frac{|\xi|}{2R}))*{\mathcal F}^{-1}\big({\overline{
 {\mathcal F}(\nabla \partial
 ^{\alpha-\beta}{\bar u}^\varepsilon\partial^\alpha
 \psi_1^\varepsilon)}}\big)\|_{L^1}
\end{eqnarray*}
Young inequality and the fact $$(*)\quad\quad\quad{\mathcal
C}^{\frac{1}{4}}\hookrightarrow L^\infty$$ we have
\begin{eqnarray*}
|\int
 \partial^\beta w_R^\varepsilon.\nabla \partial
 ^{\alpha-\beta}{\bar u}^\varepsilon\partial^\alpha
 \psi_1^\varepsilon|&\leq  & C\|\partial^\beta w_R^\varepsilon\|_{{\mathcal
C}^{s-3+\frac{1}{4}}}
 \|{\mathcal F}^{-1}(\chi(\frac{|\xi|}{2R}))\|_{L^1}\|\nabla \partial
 ^{\alpha-\beta}{\bar u}^\varepsilon\partial^\alpha
 \psi_1^\varepsilon\|_{L^2}\\
 &\leq  &C_R\|w_R^\varepsilon\|_{{\mathcal
C}^{s-3+\frac{1}{4}}}
 \|\nabla \partial
 ^{\alpha-\beta}{\bar u}^\varepsilon\|_{L^\infty}\|\partial^\alpha
 \psi_1^\varepsilon\|_{L^2}\\
 \end{eqnarray*}
the property $(*)$ imply $$ |\int
 \partial^\beta w_R^\varepsilon.\nabla \partial
 ^{\alpha-\beta}{\bar u}^\varepsilon\partial^\alpha
 \psi_1^\varepsilon|\leq C_R\|w_R^\varepsilon\|_{{\mathcal
C}^{s-3+\frac{1}{4}}}
 \|{\bar u}^\varepsilon\|_{{\mathcal
C}^{s-2+\frac{1}{4}}}\|\psi_1^\varepsilon\|_{H^{s-3}}$$ the Sobolev
injection $$\quad\quad H^{\sigma}({\mathbb R}^2
)\hookrightarrow{\mathcal C}^{\sigma-1}({\mathbb R}^2)$$ imply
\begin{eqnarray}
|\int
 \partial^\beta w_R^\varepsilon.\nabla \partial
 ^{\alpha-\beta}{\bar u}^\varepsilon\partial^\alpha
 \psi_1^\varepsilon|&\leq  &C_R\|w_R^\varepsilon\|_
 {{\mathcal C}^{s-3+\frac{1}{4}}}\|{\bar u}^\varepsilon\|_{H^{s-\frac{3}{2}+\frac{1}{4}}}
 \|\psi_1^\varepsilon\|_{H^{s-3}}\nonumber\\
&\leq  &C_R\|w_R^\varepsilon\|_{{\mathcal C}^{s-3+\frac{1}{4}}}
 \|\psi_1^\varepsilon\|_{H^{s-3}}\nonumber\\
&\leq  &C_R\|w_R^\varepsilon\|_{{\mathcal C}^{s-3+\frac{1}{4}}}^2+
 \|\psi_1^\varepsilon\|_{H^{s-3}}^2.\label{J1}\\
 \nonumber
\end{eqnarray}
$\bullet$ Study of ${\tilde J}$\\
\begin{eqnarray*}
{\tilde J}&=&|({\tilde w}_R^\varepsilon.\nabla {\bar u
}^\varepsilon;\psi_1^\varepsilon)_{H^{s-3}}|\\
&=&|\sum_{|\alpha|\leq{s-3}}\int
\partial^\alpha( {\tilde w}_R^\varepsilon.\nabla {\bar u}^\varepsilon)\partial^\alpha
\psi_1^\varepsilon|\\
&\leq
&\sum_{|\alpha|\leq{s-3}}\sum_{\beta\leq\alpha}C^\beta_\alpha|\int
\partial^\beta{\tilde w}_R^\varepsilon.\nabla \partial
^{\alpha-\beta}{\bar u}^\varepsilon\partial^\alpha
\psi_1^\varepsilon|\\
&\leq
&\sum_{|\alpha|\leq{s-3}}\sum_{\beta\leq\alpha}C^\beta_\alpha\|\nabla
\partial ^{\alpha-\beta}{\bar u}^\varepsilon\|_{L^\infty}
\|\partial^\beta{\tilde w}_R^\varepsilon\|_{L^2} \|\partial^\alpha
\psi_1^\varepsilon\|_{L^2}\\ &\leq  &C\|{\tilde
w}_R^\varepsilon\|_{H^{s-3}}
\|\psi_1^\varepsilon\|_{{H^{s-3}}}\sum_{|\alpha|\leq{s-3}}
\sum_{\beta\leq\alpha}C^\beta_\alpha\|\nabla \partial
^{\alpha-\beta}{\bar u}^\varepsilon\|_{L^\infty}\\
&\leq  &C\|{\tilde w}_R^\varepsilon\|_{H^{s-3}}
\|\psi_1^\varepsilon\|_{{H^{s-3}}}\|{\bar
u}^\varepsilon\|_{H^{s-2+\frac{1}{4}}}\\
\\
&\leq  &C\|{\tilde w}_R^\varepsilon\|_{H^{s-3}}
\|\psi_1^\varepsilon\|_{{H^{s-3}}}.\\
\end{eqnarray*}
The fact ${\mathcal F}({\tilde w}_R^\varepsilon)\equiv 0\quad
\mbox{in}\quad B(0,R)$, imply
\begin{eqnarray*}
\|{\tilde
w}_R^\varepsilon\|_{H^{s-3}}&\leq  &R^{-3}\|w^\varepsilon\|_{H^{s}}\\
&\leq  &C. R^{-3}\\
\end{eqnarray*}
then
\begin{eqnarray}{\tilde J}&\leq  &C.
R^{-3}\|\psi_1^\varepsilon\|_{{H^{s-3}}}\nonumber\\ &\leq  &C.
R^{-2}\|\psi_1^\varepsilon\|_{{H^{s-3}}}\nonumber\\
&\leq  &C.R^{-4}+\|\psi_1^\varepsilon\|_{{H^{s-3}}}^2.\label{J2}\\
\nonumber\end{eqnarray} Using (\ref{ef})...(\ref{J2}) we obtain
$$\|\psi^\varepsilon(t)\|^2_{H^{s-3}}\leq
C.R^{-4}+C_R\int_0^{T_1}\|w_R^\varepsilon(\tau)\|_{{\mathcal
C}^{s-3+\frac{1}{4}}}^2d\tau+C\int_0^t\|\psi^\varepsilon(\tau)\|_{H^{s-3}}d\tau,
$$ Gronwall lemma imply $$\|\psi^\varepsilon(t)\|^2_{H^{s-3}}\leq
\Big(C.R^{-4}+C_R\int_0^{T_1}\|w_R^\varepsilon(\tau)\|_{{\mathcal
C}^{s-\frac{11}{4}}}^2d\tau\Big)\exp(CT_1),$$
 The proofs of Theorem \ref{t45} and \ref{t47} imply that
$$w^\varepsilon_R\rightarrow 0\quad in\quad
L^4([0,T_1],{\mathcal C}^{s'-{\frac{3}{2}}}),\quad \forall s'< s;$$
we can deduce the desired result.\endproof
\section{ Appendix}
\noindent In this section we shall prove the product laws stated in
Lemmas \ref{l41}, \ref{l47} and \ref{l48}.
\subsection{Proof of Lemma \ref{l41}}
We suppose that $a,b, c$ are three vectors fildes in
$H^\infty({\mathbb R}^3)$ such that $\mbox{div}~a=0$.\\\\  For the
first point:\\ We write $$|{< a.\nabla b, b>
}_{H^\sigma(\mathbb{R}^3)}|=|\sum_{|\alpha|\leq\sigma}\int\partial^{\alpha}(
a.\nabla b)\partial^{\alpha}b|,$$ then
\begin{eqnarray*}
|{< a.\nabla b, b> }_{H^\sigma(\mathbb{R}^3)}|&\leq
&\sum_{|\alpha|\leq\sigma}|\int\partial^{\alpha}(a.\nabla
b)\partial^{\alpha}b|\\&\leq
&\sum_{|\alpha|\leq\sigma}\sum_{\beta\leq\alpha}C_{\alpha}^{\beta}
A_{\alpha,\beta},
\end{eqnarray*}
where $$ A_{\alpha,\beta}:=|\int(\partial^{\beta}a.\nabla
\partial^{\alpha-\beta}b)\partial^{\alpha}b|.$$
The most important term is for $|\alpha|=\sigma$.\\ $\bullet$~
$\mbox{div}~a=0$ imply that $A_{\alpha,\alpha}=0$.\\$\bullet$ For
$0< |\beta|< \sigma$, we are $$
\begin{array}{cllll}
A_{\alpha,\alpha}&\leq  &\displaystyle\int
|\partial^{\alpha-\beta}a||\nabla\partial^{\beta}b||\partial^{\alpha}b|\\
&\leq
&\|\partial^{\alpha-\beta}a\|_{L^\infty}\|\nabla\partial^{\beta}b\|_{L^2}
\|\partial^{\alpha}b\|_{L^2}\\
&\leq  &\|\partial^{\alpha-\beta}a\|_{{\mathcal
C}^{\frac{1}{4}}}\|\nabla b\|_{H^{\sigma-1}}^2\\ &\leq  & \|\nabla
a\|_{{\mathcal C}^{\sigma-\frac{11}{4}}}\|\nabla
b\|_{H^{\sigma-1}}^2\\&\leq  & C\|\nabla a\|_{H^{\sigma-1}}\|\nabla
b\|_{H^{\sigma-1}}^2\\
\end{array}
$$
\\
$\bullet$ For $\beta=0$, $$
\begin{array}{cllll}
A_{\alpha,\beta}&\leq  &\displaystyle\int
|\partial^{\alpha}a||\nabla b||\partial^{\alpha}b|\\
&\leq  &\|\partial^{\alpha}a\|_{L^2}\|\nabla b\|_{L^\infty}
\|\partial^{\alpha}b\|_{L^2}\\ &\leq  &\|\nabla
a\|_{H^{\sigma-1}}\|\nabla b\|_{{\mathcal C}^{\frac{1}{4}}}\|\nabla
b\|_{H^{\sigma-1}}\\ &\leq &C\|\nabla a\|_{H^{\sigma-1}} \|\nabla
b\|_{H^{\frac{7}{4}}}\|\nabla b\|_{H^{\sigma-1}}\\&\leq  & C\|\nabla
a\|_{H^{\sigma-1}}\|\nabla
b\|_{H^{\sigma-1}}^2\\
\end{array}
$$
For the second point, we write $$\begin{array}{cllll} |{< a\nabla
b,c>}_{H^\sigma}+{< a\nabla
c,b>}_{H^\sigma}|&=&\displaystyle\sum_{\alpha\in{\mathbb N }^3\atop
|\alpha|\leq \sigma}\big(\int\partial^\alpha (a\nabla
b)\partial^\alpha c
+\int\partial^\alpha (a\nabla c)\partial^\alpha b\big)\\
&\leq&\displaystyle\sum_{\alpha\in{\mathbb N }^3\atop |\alpha|\leq
\sigma}\sum_{\beta\leq\alpha}C^\beta_\alpha B_{\alpha,\beta},\\
\end{array}$$
where
$$B_{\alpha,\beta}:=|\int(\partial^{\alpha-\beta}a.\nabla
\partial^{\beta} b)
\partial^\alpha c +\int(\partial^{\alpha-\beta}a.\nabla
\partial^{\beta} c)\partial^\alpha
b|$$

 The most important term is for $|\alpha|=\sigma$.\\
$\bullet$~$\mbox{div}~a=0$ imply that $B_{\alpha,\alpha}=0$.\\
$\bullet$ For $\beta\neq\alpha$ we apply the first step.\endproof
\subsection{Proof of Lemma \ref{l47}}
We write $$\begin{array}{cllll}|{< a.\nabla b,
c>}_{H^\sigma(\mathbb{R}^3)}|&=&
|\displaystyle\sum_{|\alpha|\leq\sigma}\int\partial^{\alpha}(a.\nabla
b)\partial^{\alpha}c|\\
&=&|\displaystyle\sum_{|\alpha|\leq\sigma}\int\partial^{\alpha}(a.\nabla
b)\partial^{\alpha}c|\\
&=&|\displaystyle\sum_{|\alpha|\leq\sigma}\sum_{\beta\leq\alpha}
C^\beta_\alpha\int(\partial^{\alpha-\beta} a.\nabla
\partial^{\beta}b)\partial^{\alpha}c|\\
&\leq
&\displaystyle\sum_{|\alpha|\leq\sigma}\sum_{\beta\leq\alpha}C^\beta_\alpha
D_{\alpha,\beta},\\
\end{array}$$
where
$$D_{\alpha,\beta}:=|\displaystyle\int(\partial^{\alpha-\beta}
a.\nabla
\partial^{\beta}b)\partial^{\alpha}c|$$
the Cauchy-Schwarz inequality imply
$$
\begin{array}{cllll}
D_{\alpha,\beta}&\leq  &\|\partial^{\alpha-\beta} a.\nabla
\partial^{\beta}b\|_{L^2}\|\partial^{\alpha}c\|_{L^2}\\
&\leq
&\|\partial^{\alpha}c\|_{L^2}\min\big(\|\partial^{\alpha-\beta}
a\|_{L^2}\|\nabla
\partial^{\beta}b\|_{L^\infty},\|\partial^{\alpha-\beta}
a\|_{L^\infty}\|\nabla
\partial^{\beta}b\|_{L^2}\big)\\
&\leq
&\|\partial^{\alpha}c\|_{L^2}\min\big(\|\partial^{\alpha-\beta}
a\|_{L^2}\|\nabla
\partial^{\beta}b\|_{{\mathcal C}^\frac{1}{4}},\|\partial^{\alpha-\beta}
a\|_{{\mathcal C}^\frac{1}{4}}\|\nabla
\partial^{\beta}b\|_{L^2}\big)\\
&\leq  &\|c\|_{H^\sigma}\min\big(\| a\|_{H^\sigma}\|\nabla
b\|_{{\mathcal C}^{\sigma+\frac{1}{4}}},\|a\|_{{\mathcal
C}^{\sigma+\frac{1}{4}}}\|\nabla b\|_{H^\sigma}\big)\\
\end{array}
$$ This achieves the proof.\endproof
\subsection{Proof of Lemma \ref{l48}}
We write $$\begin{array}{cllll}|{<{\bar a}.\nabla b,
c>}_{H^\sigma(\mathbb{R}^3)}+{<{\bar a}.\nabla c,
b>}_{H^\sigma(\mathbb{R}^3)}|&=&
|\displaystyle\sum_{|\alpha|\leq\sigma}\int\partial^{\alpha}({\bar
a}.\nabla b)\partial^{\alpha}c+\int\partial^{\alpha}({\bar
a}.\nabla c)\partial^{\alpha}b|\\
&=&|\displaystyle\sum_{|\alpha|\leq\sigma}\sum_{\beta\leq\alpha}
C^\beta_\alpha\Big(\int(\partial^{\alpha-\beta}{\bar a}.\nabla
\partial^{\beta}b)\partial^{\alpha}c\\&&+\displaystyle\int(\partial^{\alpha-\beta}{\bar a}.\nabla
\partial^{\beta}c)\partial^{\alpha}b\Big)|\\
&\leq
&\displaystyle\sum_{|\alpha|\leq\sigma}\sum_{\beta\leq\alpha}C^\beta_\alpha
E_{\alpha,\beta},\\
\end{array}$$
where
$$E_{\alpha,\beta}:=|\displaystyle\int(\partial^{\alpha-\beta}{\bar
a}.\nabla
\partial^{\beta}b)\partial^{\alpha}c+\int(\partial^{\alpha-\beta}{\bar a}.\nabla
\partial^{\beta}c)\partial^{\alpha}b|.$$
The most important term is for $|\alpha|=\sigma$.\\
$\bullet$ $\mbox{div}_h~{\bar a}=0$ imply that $E_{\alpha,\alpha}=0$.\\
For $\beta\neq\alpha $ we are $E_{\alpha,\beta}
\leq E_{\alpha,\beta}^{(1)}+E_{\alpha,\beta}^{(2)},$\\
where
$E_{\alpha,\beta}^{(1)}=|\displaystyle\int(\partial^{\alpha-\beta}{\bar
a}.\nabla
\partial^{\beta}b)\partial^{\alpha}c|$ and $E_{\alpha,\beta}^{(2)}=|\displaystyle\int(\partial^{\alpha-\beta}{\bar a}.\nabla
\partial^{\beta}c)\partial^{\alpha}b|$.\\
 $\bullet$
$0\leq |\beta|\leq\sigma-2$.\\In this case we have $$
E_{\alpha,\beta}^{(1)}\leq \|\partial^{\beta}{\overline
a}\|_{L^\infty(\mathbb{R}^2)}\|\nabla
\partial^{\alpha-\beta}b\|_{L^2(\mathbb{R}^3)}\|\partial^{\alpha}c\|_{L^2(\mathbb{R}^3)}.$$
Using the injections $$
H^{1+\frac{1}{4}}(\mathbb{R}^2)\hookrightarrow{\mathcal
C}^{\frac{1}{4}}(\mathbb{R}^2) \hookrightarrow
L^\infty(\mathbb{R}^2)$$ we deduce $$E_{\alpha,\beta}^{(1)}\leq
C\|{\bar a}\|_{H^\sigma(\mathbb{R}^2)}\|\nabla b
\|_{H^{\sigma-1}(\mathbb{R}^3)}\|\nabla
c\|_{H^{\sigma-1}(\mathbb{R}^3)}.$$ $\bullet$ $\sigma-1\leq
|\beta|\leq\sigma.$\\ If we denote $x=(x_1,x_2,x_3):=(x_h,x_3)$,
then
\begin{eqnarray*}
E_{\alpha,\beta}^{(1)}&\leq
&\int_{x_3}\Big(\int_{x_h}|\partial^{\beta}{\overline
a}(x_h)||\nabla
\partial^{\alpha-\beta}b(x_h,x_3)||\partial^{\alpha}c(x_h,x_3)|\Big)\\
&\leq  & \|\partial^{\beta}{\overline
a}\|_{L^2(\mathbb{R}^2)}\Big(\int_{x_3}\|\nabla
\partial^{\alpha-\beta}b(x_3)\|_{L^\infty(\mathbb{R}^2)}^2dx_3
\Big)^{\frac{1}{2}}\Big(\int_{x_3}
\|\partial^{\alpha}c(x_3)\|_{L^2(\mathbb{R}^2)}^2dx_3\Big)^{\frac{1}{2}}.
\end{eqnarray*}
But
\begin{eqnarray*}
\|\nabla
\partial^{\alpha-\beta}b(x_3)\|_{L^\infty(\mathbb{R}^2)}&\leq  &C\|\nabla
\partial^{\alpha-\beta}b(x_3)\|_{H^{\frac{3}{2}}(\mathbb{R}^2)}\\
&\leq  &C \|\nabla
\partial^{\alpha-\beta}b(x_3)\|_{H^2(\mathbb{R}^2)},
\end{eqnarray*}
$$ \Big(\int_{x_3}
\|\partial^{\alpha}c(x_3)\|_{L^2(\mathbb{R}^2)}^2dx_3\Big)^{\frac{1}{2}}=
\|\partial^{\alpha}c\|_{L^2(\mathbb{R}^3)},$$ and $$
\Big(\int_{x_3}\|\nabla
\partial^{\alpha-\beta}b(x_3)\|_{H^2(\mathbb{R}^2)}^2\Big)^{\frac{1}{2}}\leq
\|\nabla
\partial^{\alpha-\beta}b\|_{H^2(\mathbb{R}^3)}
$$ so $$ E_{\alpha,\beta}\leq \|{\bar
a}\|_{H^\sigma(\mathbb{R}^2)}\|\nabla
b\|_{H^{\sigma-1}(\mathbb{R}^3)}\|c\|_{H^\sigma(\mathbb{R}^3)}.$$
For the second estimate we remark that, if $c=b$,
$$\mbox{div}_h~{\bar a }=0\quad \mbox{imply}\quad
E_{\alpha,\alpha}=0\quad\mbox{for all} \quad |\alpha|\leq \sigma.$$
For the third point we write $$|{< b.\nabla {\overline a},
c>}_{H^\sigma(\mathbb{R}^3)}|=|\sum_{|\alpha|\leq\sigma}\int\partial^{\alpha}(b.\nabla
{\overline a})\partial^{\alpha}c|,$$ then $$ |{< b.\nabla {\overline
a},
c>}_{H^\sigma(\mathbb{R}^3)}|\leq\sum_{|\alpha|\leq\sigma}\sum_{\beta\leq\alpha}
C_{\alpha}^{\beta} F_{\alpha,\beta}, $$ where $$
F_{\alpha,\beta}:=|\int(\partial^{\beta}b.\nabla
\partial^{\alpha-\beta}{\overline a})\partial^{\alpha}c|.$$
The most important case is when $|\alpha|=\sigma$.\\
$\bullet$ For $3\leq |\beta|\leq\sigma$, we can write
$$F_{\alpha,\beta}\leq\|\nabla\partial^{\alpha-\beta}{\overline
a}\|_{L^\infty(\mathbb{R}^2)}\|\partial^{\beta}b\|_{L^2(\mathbb{R}^3)}
\|\partial^{\alpha}c\|_{L^2(\mathbb{R}^3)}.$$ But
$$\|\nabla\partial^{\alpha-\beta}{\overline
a}\|_{L^\infty(\mathbb{R}^2)}\leq C\|{\overline
a}\|_{H^\sigma(\mathbb{R}^2)}.$$ So $$F_{\alpha,\beta}\leq
C\|{\overline
a}\|_{H^\sigma(\mathbb{R}^2)}\|b\|_{H^\sigma(\mathbb{R}^3)}\|c\|_{H^\sigma(\mathbb{R}^3)}.$$
$\bullet$ For $|\beta|\leq 2$, we have by the Cauchy-Schwarz
inequality
\begin{eqnarray*}
F_{\alpha,\beta}&\leq
&\int_{x_3}\Big(\int_{x_h}|\nabla\partial^{\alpha-\beta}{\overline
a}(x_h)||\partial^{\beta}b(x_h,x_3)||\partial^{\alpha}c(x_h,x_3)|\Big)\\
&\leq  & \|\nabla\partial^{\alpha-\beta}{\overline
a}\|_{L^2(\mathbb{R}^2)}\Big(\int_{x_3}\|\partial^{\beta}b(x_3)\|_{L^\infty(\mathbb{R}^2)}^2dx_3
\Big)^{\frac{1}{2}}\Big(\int_{x_3}
\|\partial^{\alpha}c(x_3)\|_{L^2(\mathbb{R}^2)}^2dx_3\Big)^{\frac{1}{2}}\\
&\leq  &C\|{\overline
a}\|_{H^{\sigma+1}(\mathbb{R}^2)}\|b\|_{H^\sigma(\mathbb{R}^3)}
\|c\|_{H^\sigma(\mathbb{R}^3)}.
\end{eqnarray*}
This end the proof of lemma \ref{l48}.\endproof

\end{document}